          \documentclass{cmslatex}
\usepackage[paperwidth=7in, paperheight=10in, margin=.875in]{geometry}

\usepackage{graphicx}
\usepackage{amssymb}
\usepackage{amsmath}
\usepackage{epstopdf}
\usepackage{subcaption}
\usepackage{tikz}
\usepackage{multirow}
\newcommand{\Rbb}{\mathbb{R}}

\newcommand{\bm}[1]{\mbox{\boldmath${#1}$}}

\newcommand{\domain}{\Omega}
\newcommand{\cdomain}{\bar{\Omega}}
\newcommand{\boundary}{\partial \domain}

\newcommand{\U}{\mathcal{U}}

\newcommand{\x}{\boldsymbol{x}}

\newcommand{\bn}{\bm{n}}
\newcommand{\bp}{\bm{p}}
\newcommand{\bq}{\bm{q}}
\newcommand{\bu}{\boldsymbol{u}}
\newcommand{\bv}{\bm{v}}
\newcommand{\br}{\boldsymbol{r}}
\newcommand{\bV}{\bm{V}}
\newcommand{\bw}{\bm{w}}
\newcommand{\w}{\bm{w}}
\newcommand{\y}{\bm{y}}

\newcommand{\rhoA}{\rho^{\textsc{a}}}
\newcommand{\rhoB}{\rho^{\textsc{b}}}
\newcommand{\buA}{\bm{u}^{\textsc{a}}}
\newcommand{\buB}{\bm{u}^{\textsc{b}}}
\newcommand{\bvA}{\bm{v}^{\textsc{a}}}

\newcommand{\vA}{v^{\textsc{a}}}
\newcommand{\vB}{v^{\textsc{b}}}
\newcommand{\varphiA}{\varphi^\textsc{a}}
\newcommand{\varphiB}{\varphi^\textsc{b}}
\newcommand{\PhiA}{\Phi^\textsc{a}}
\newcommand{\PhiB}{\Phi^\textsc{b}}
\newcommand{\JA}{J^\textsc{a}}
\newcommand{\JB}{J^\textsc{b}}
\newcommand{\fA}{f^\textsc{a}}
\newcommand{\fB}{f^\textsc{b}}

\newcommand{\dfadpsi}{\frac{d\fA}{d\psi}}
\newcommand{\dfaddpsi}{\frac{d^2\fA}{d\psi^2}}
\newcommand{\rhobar}{\bar{\rho}}
\newcommand{\V}{\mathcal{V}}

\newcommand{\argmax}{\mathop{\mathrm{argmax}}}

\renewcommand{\theequation}{\arabic{section}.\arabic{equation}}

\begin{document}
          \title{Anisotropic Challenges in Pedestrian Flow Modeling\thanks{Both authors' work was supported in part by the National Science Foundation grants DMS-1016150 and DMS-1738010.  The first author's work was also supported 
          by the National Science Foundation grants DMS-0739164 and DMS-1645643.}}

          \author{Elliot Cartee\thanks{Department of Mathematics, Cornell University, Ithaca, NY 14853 (evc34@cornell.edu)}
          \and{Alexander Vladimirsky\thanks{Department of Mathematics, Cornell University, Ithaca, NY 14853 (vladimirsky@cornell.edu)}}}

         \pagestyle{myheadings} \markboth{Anisotropic Crowd Modeling}{E. Cartee and A. Vladimirsky} \maketitle

          \begin{abstract}
Macroscopic models of crowd flow incorporating individual pedestrian choices present many analytic and computational challenges.
Anisotropic interactions are particularly subtle, both in terms of describing the correct ``optimal'' direction field for the pedestrians
and ensuring that this field is uniquely defined.  We develop sufficient conditions, which establish a range of ``safe'' densities and
parameter values for each model.  We illustrate our approach by analyzing several established intra-crowd and inter-crowd models.
For the two-crowd case, we also develop sufficient conditions for the uniqueness of Nash Equilibria in the resulting non-zero-sum game.
          \end{abstract}
\begin{keywords}
pedestrian dynamics, conservation laws, Hamilton-Jacobi equations, Hughes' model,
anisotropic speed profiles, minimum time problem, non-zero-sum differential games, Nash equilibrium
\end{keywords}

 \begin{AMS}
35L65, 49N90, 91D10, 65N22.
\end{AMS}

\section{Introduction.}
\label{s:Intro}
Modeling the dynamics of crowds is a very active area in traffic engineering.
A good overview and classification of most popular approaches can be found in \cite{Bellomo_SIAMRev}.
For some of the macroscopic models incorporating individual choices of pedestrians, significant efforts have already been invested to build efficient implementations and conduct extensive numerical experiments (e.g., \cite{Huang_2009, Jiang_2009, xiong2011, xiong2011lanes, jiang2012numerical, jiang2015higher, Cristiani_2015,  Hartmann_2014, hanseler2015, jiang2016macroscopic, jiang2016comparison, carrillo2016improved}), despite the fact that some of the theoretical underpinnings are still missing.  This includes both the well-posedness of PDE systems and the convergence under grid refinement of approximate solutions.  Some practitioners may even argue that the issues of numerical convergence are not particularly relevant \cite{Hartmann_2014} since these PDEs are only approximations of crowd dynamics valid in intermediate asymptotics and the actual traffic engineering applications assume a specific length-scale.  In the absence of rigorous error estimates,
the models are judged instead phenomenologically -- based on the perceived realism of their numerical predictions and on their ability to generate some of the self-organized behavior (e.g., ``lane formation'') observed in real crowds.  This tendency seems 
problematic given the blurring of boundaries between qualitative models (used to highlight the underlying mechanisms) and quantitative models (used to make specific predictions or even change the shape of the environment -- to avoid stampedes or to improve the chances of a successful evacuation).

In macroscopic models the crowd is represented as a density function $\rho(\x,t)$, whose evolution is described by a conservation law.
A typical (hyperbolic) version is
\begin{equation}
\rho_t(\x,t) \; + \; \nabla \cdot (\rho(\x,t)\bV(\x,t)) \; = \; 0, \qquad \qquad \x \in \Omega \subseteq \Rbb^2, \, t>0,
\label{Hyperbolic1}
\end{equation}
where $\bV$ is the velocity of motion for all pedestrians located at $\x \in \Omega$ at the time $t.$
The initial density $\rho(\x,0)$ is assumed to be known and additional boundary conditions are used to specify the obstacles and inflow/outflow rates through the exits.  In many approaches, the velocity field $\bV$ is obtained by modeling pedestrians' individual choices -- they plan their paths 
taking into account the current (and possibly future) crowd density on
different parts of
$\Omega$; see, for example, \cite{Huang_2009, Jiang_2009, Cristiani_2015}.
This path optimization is performed by solving a nonlinear Hamilton-Jacobi equation coupled to \eqref{Hyperbolic1}.  The well-posedness of the resulting PDE systems is far from trivial and is an area of active research \cite{DIFRANCESCO2011, AMADORI2012, Priuli_MFG}.
Our goals in this paper are much more modest: we simply aim to verify the internal/logical model consistency, by which we mean that the velocity field
$\bV$ should be uniquely defined almost everywhere on $\Omega$.

The internal consistency requirement is trivially satisfied
when the pedestrians' velocities are isotropic; i.e., they can move in any direction with the same speed dependent on $\rho$ at the current position $(\x,t$).  This is precisely the setting described by Hughes \cite{Hughes_2002} and reviewed in section \ref{s:OptimalControl}.
However, anisotropic velocities arise naturally in many generalizations: due to the landscape (e.g., walking on a steep hill), or interactions of multiple crowds (e.g., \cite{Jiang_2009, jiang2012numerical, jiang2015higher, hanseler2015}), or non-local dependencies on $\rho$ (e.g., \cite{Cristiani_2015}).
Regardless of the source of anisotropy, the notion of ``optimal direction of motion'' becomes far more delicate and can result in subtle mistakes in computational implementations.
A {\em velocity profile} $\V(\x,t)$ represents the set of all velocities that a pedestrian can use at a specific time and location; its formal definition 
is included in section \ref{s:SpeedProfiles}, where we show that its geometric properties determine the optimal pedestrian choices.
We also show that the strict convexity of  $\V$ is needed to guarantee that anisotropic models are internally consistent.

Unfortunately, computational experiments by themselves are insufficient to uncover model inconsistencies since they might be masked by the implementation details.  Throughout this paper, we conduct a 
review of implicit assumptions in prior models. 
However, our real focus is on developing minimal consistency conditions, which should be easy to verify -- either beforehand or in conjunction with any numerical implementation.  Whether or not our conditions are actually satisfied in existing models depends on the crowd densities and parameter values.  This brings us to a discussion of some of the experimental work to determine the latter (see section \ref{ss:SpeedProfiles_Inter}).

With two interacting crowds, we argue that models are consistent only if there is no need to specify which crowd ``goes first'' in choosing its velocity field.
This is essentially the requirement for the uniqueness of a Nash equilibrium in a non-zero-sum differential game,
and in section \ref{s:Nash} we derive the sufficient conditions to guarantee it.

In section \ref{s:NumericalImplementation} we discuss the challenges of building efficient numerical algorithms for anisotropic models.
For a single crowd with non-local interactions, we propose a new approach to decrease the computational cost by using a linear approximation of crowd density.
For two-crowd models, we discuss the use of iterative solvers for a coupled system of Hamilton-Jacobi-Isaacs PDEs and include the results of several computational experiments.
We conclude by discussing several desirable extensions in section \ref {s:Conclude}.

\section{Optimal Control Formulation.}
\label{s:OptimalControl}
Finding time-optimal trajectories to a target is among the most studied problems in optimal control theory.
Here we provide a brief overview of the dynamic programming approach, focusing on the geometric interpretation and referring readers to
 \cite{bardicapuzzodolcetta, bressan2007introduction} for a detailed treatment.
Suppose the local velocity depends on one's location and the chosen control value; i.e., $\bv: \cdomain \times \U \mapsto \Rbb^2$,
where $\U$ is a suitable compact set of control values.  In our context, it will be enough to identify controls as the pedestrians' preferred directions of motion; i.e.,
$\U = S^1 = \{ \bu \in \Rbb^2 \, \mid \, | \bu | = 1 \}.$
We will focus on two types of pedestrian velocities:
\begin{align}
\label{eq:geom_velocity}
\bv(\x,\bu) \; &= \; v(\x,\bu) \bu,  \; \\
\nonumber
& \text{ where  $v>0$ is the speed of motion through $\x$ in the direction $\bu$; and } \\
\label{eq:additive_velocity}
\bv(\x,\bu) \; &= \; \bu \, + \, \bw_i(\x,\bu),  \; \\
\nonumber
& \text{ where  $\bw_i$ is a known ``interactional velocity'' perturbation. }
\end{align}

The dynamics are governed by $\y'(t) = \bv(\y(t),\bu)$, and we want to choose a feedback control function $\bu = \bu(\y(t))$ to minimize the time of travel to the target $\Gamma_d$.
The {\em value function} $\varphi(\x)$ is defined to be the minimum time to target from the starting position $\y(0) = \x \in \Omega$,
and it can be recovered as the unique viscosity solution \cite{crandall1983viscosity} of the Hamilton-Jacobi-Bellman PDE
\begin{align}
\label{eq:HJB_general}
\max_{\bu \in \U}\{-\nabla \varphi(\x) \cdot \bv(\x, \bu) \} =& \, 1, &\x \in \Omega \\
\varphi(\x) =& \, 0, &\x \in \Gamma_d \nonumber
\end{align}
with additional boundary conditions (e.g., $\varphi=+\infty$ on $\boundary \backslash \Gamma_d$) to ensure that all trajectories stay inside $\Omega$.
The characteristics of this PDE are in fact the time-optimal trajectories and a time-optimal feedback control is obtained by choosing $\bu(\x)$ from the argmax set of equation \eqref{eq:HJB_general}.



The above describes path-planning for a pedestrian whose velocity is only affected by the current location and the chosen direction.
But the dependence on crowd density can be treated similarly; e.g., with the speed $v(\rho,\bu)$ monotone decreasing in $\rho$.
If the rest of the crowd is not moving and the current density $\hat{\rho}$ is known, our pedestrian can simply use  $v(\x,\bu) = v(\hat{\rho}(\x),\bu)$
to plan her path to $\Gamma_d$.
We will refer to this approach as ``stroboscopic'' since in reality all pedestrians are moving using the same logic, but our planning is based on the current snapshot of the crowd density instead of trying to predict the correct/updated density at each point on the trajectory\footnote{The latter approach is also present in the literature \cite{hoogendoorn2003simulation,hoogendoorn2004dynamic,lachapelle2011mfg,burger2013mfg,Cristiani_2015,Priuli_MFG} and is related to the {\em Mean Field Games} \cite{lasry2006jeux,lasry2006jeux2,HuangMalhameCaines,gueant2011mfg,Gomes2014}. Many of the issues we discuss are relevant in that context as well, but we focus on the stroboscopic models for the sake of simplicity.}. So, the current $\varphi$ is only good enough to choose the initial direction of motion and has to be continuously recomputed as the density changes. Putting this all together, our system of PDEs is
\begin{subequations} \label{eq:strobo_full}
\begin{align}
\rho_t(\x,t) \; + \; \nabla \cdot (\rho(\x,t)\bV(\x,t)) &= \; 0,  \tag{\ref{eq:strobo_full}a} \\
\max\limits_{\bu \in \U} \left\{-\nabla \varphi(\x,t) \cdot \bv(\rho(\x,t), \bu) \right\} &= \, 1, \tag{\ref{eq:strobo_full}b}\\
\bu_*(\x,t) \, \in \, \argmax\limits_{\bu \in \U}
\left\{-\nabla \varphi(\x,t) \cdot \bv(\rho(\x,t), \bu) \right\}, & \tag{\ref{eq:strobo_full}c} \\
\bV(\x,t) = \bv(\rho(\x,t), \bu^*(\x,t)), & \tag{\ref{eq:strobo_full}d}
\end{align}
\end{subequations}
with the specified initial density $\rho(\x,0) = \rho_0(\x)$ and $\varphi(\x,t) = 0$ on  $\Gamma_d \times R$.

\begin{figure}[h]
\centerline{
$
\arraycolsep=-3pt
\begin{array}{ccccc}	
\includegraphics[width=.095\textwidth]{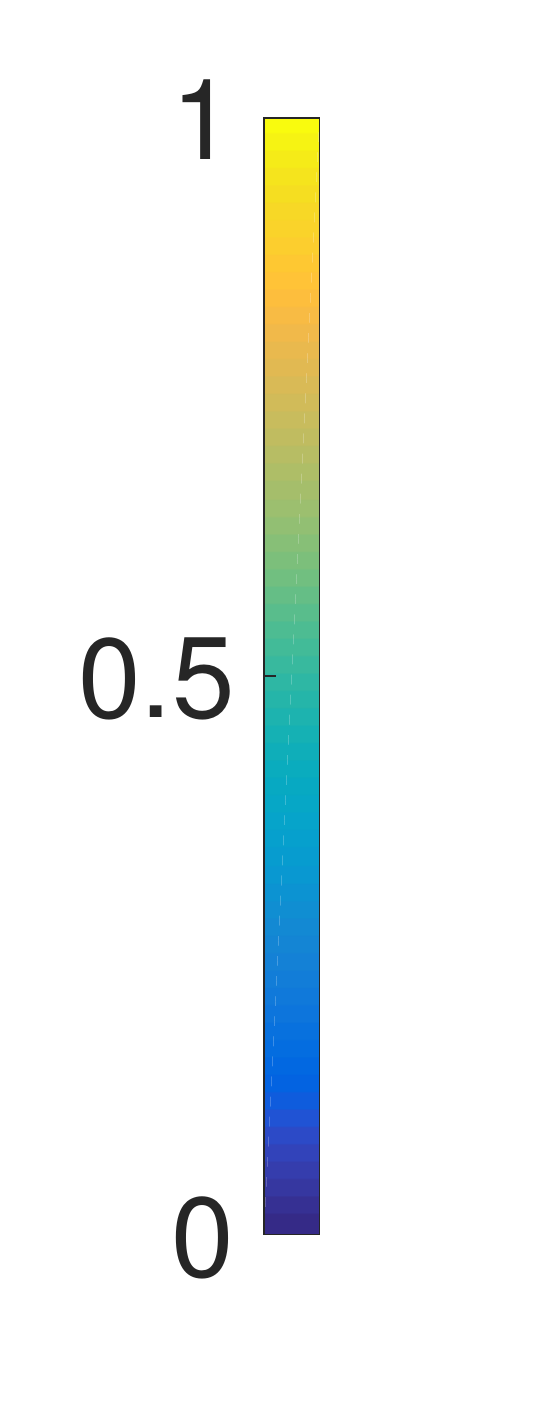} &
\includegraphics[width=.28\textwidth]{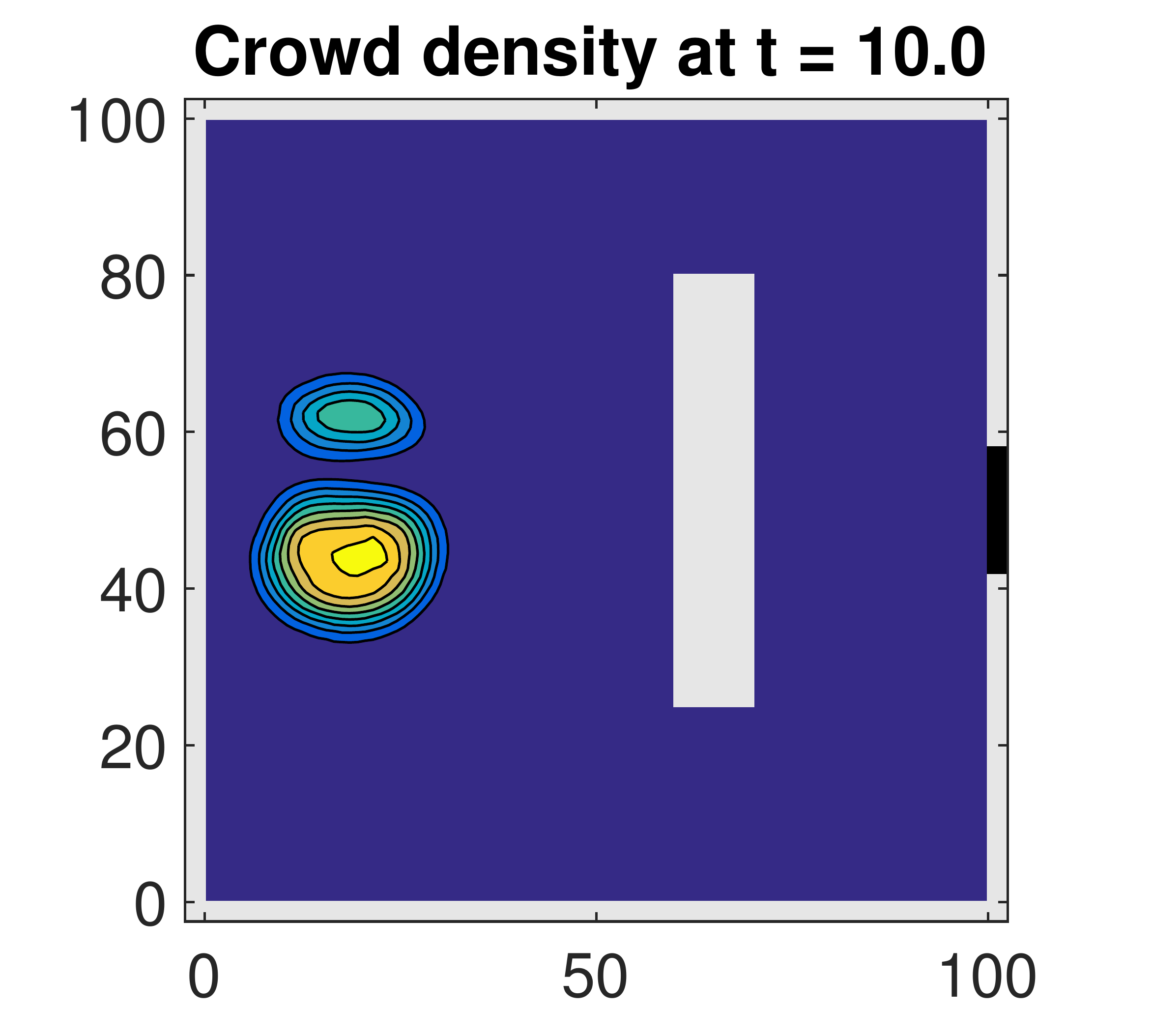} &
\includegraphics[width=.28\textwidth]{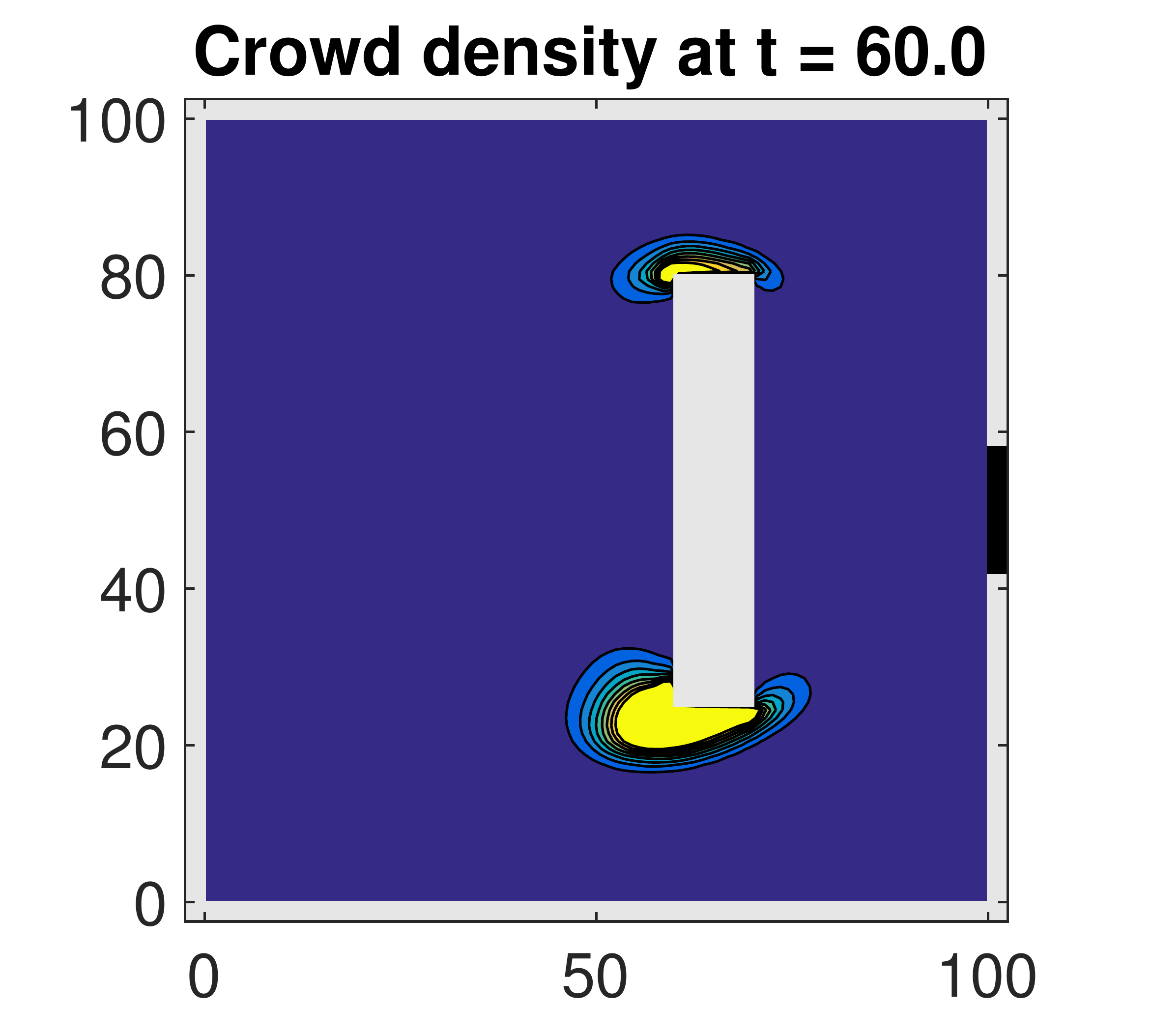} &
\includegraphics[width=.28\textwidth]{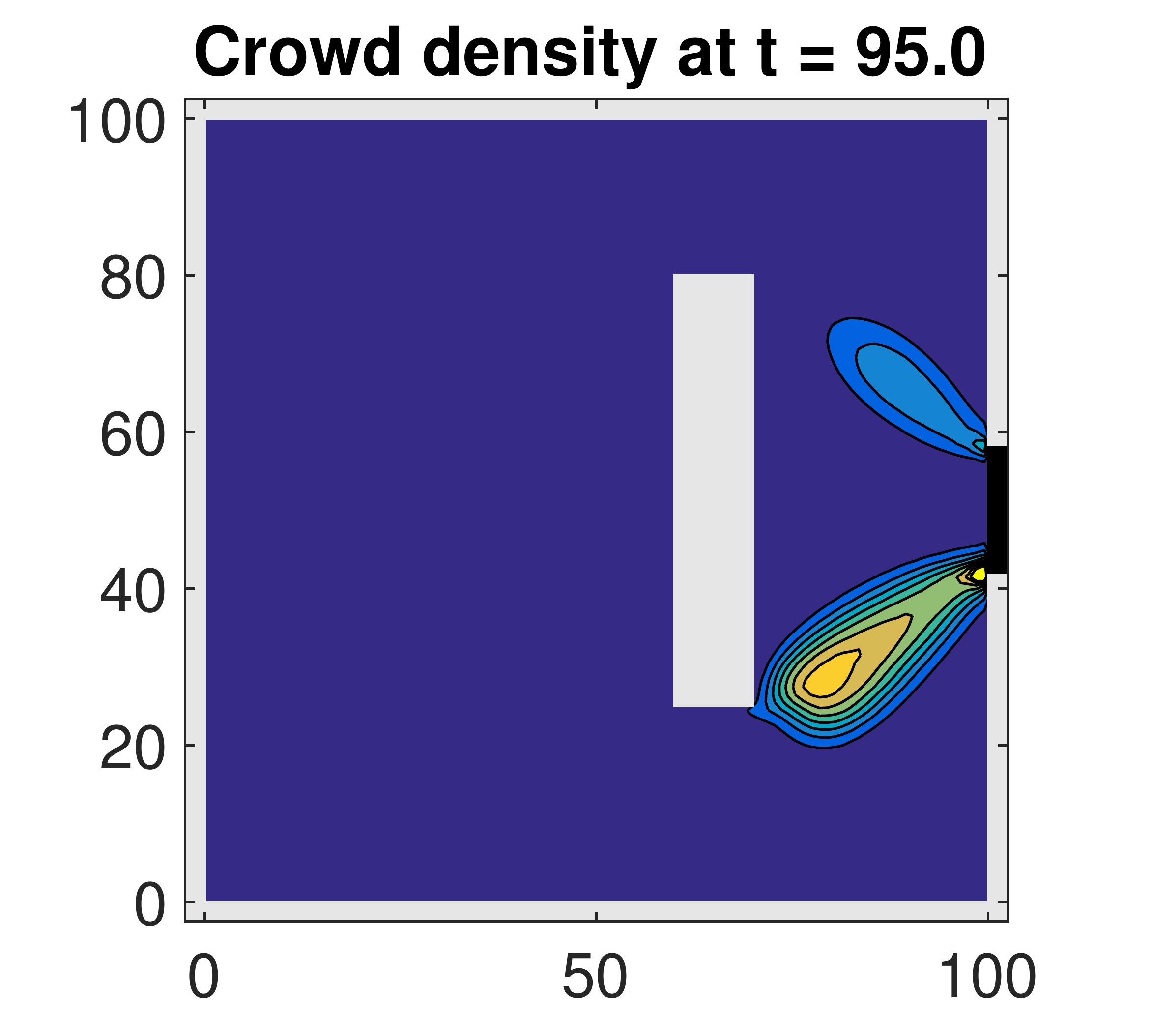} &
\includegraphics[width=.28\textwidth]{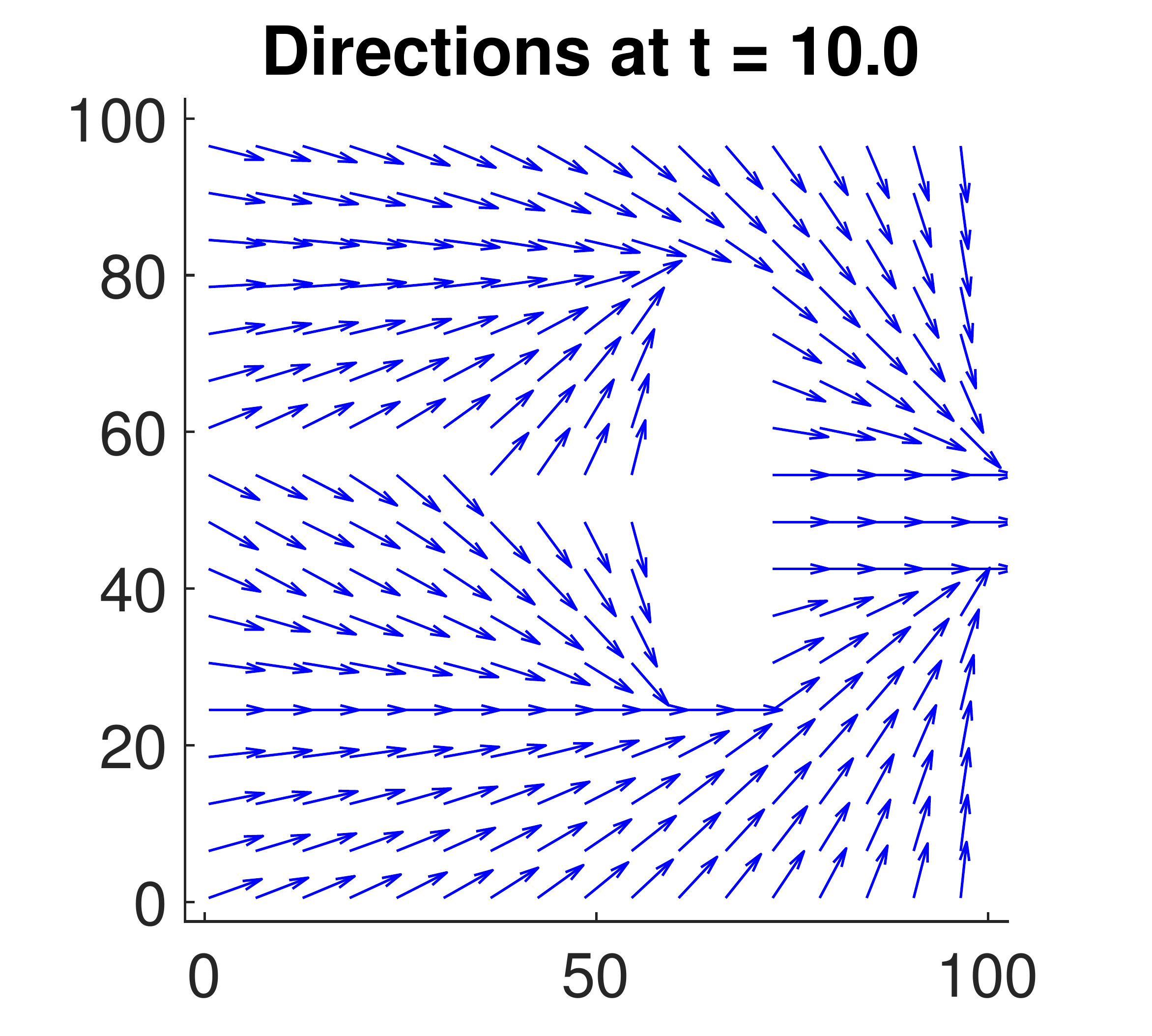} \\
& (a) & (b) & (c) & (d)
\end{array}
$
}
\caption{Simulation of a single (isotropic) crowd evacuating through a doorway (shown in black) on the right edge of the domain. Contour plots of $\rho$ for various time-slices (in a,b,c) and a quiver plot of $\bu_*$ in d. The speed function is $v(\rho)=e^{-0.075\rho^2}.$
Computed on a 100 x 100 grid with a timestep $\Delta t= 0.5$.}
\label{fig:SingleCrowd}
\end{figure}

In the isotropic case (where $v(\rho,\bu)=v(\rho)$), each $\varphi(\cdot,t)$ is a solution of the Eikonal equation:
\begin{equation}\label{eq:Eikonal}
|\nabla \varphi(\x,t)|v(\rho(\x,t)) = 1
\end{equation}
 and  the crowd velocity is
$\bV(\x,t) = - v(\rho(\x,t)) \frac{\nabla \varphi(\x,t)}{|\nabla \varphi(\x,t)|} = - v^2(\rho(\x,t)) \nabla \varphi(\x,t).$
(A typical isotropic evolution of crowd density on a domain with a single obstacle is shown in Figure~\ref{fig:SingleCrowd}.)
This is essentially the set up considered in the influential paper by Hughes \cite{Hughes_2002}, who was among the first to model
the influence of pedestrian's choices on the macroscopic dynamics\footnote{The equation introduced by Hughes was slightly more general than the time-optimal path planning -- it allowed for an additional high-density-discomfort/penalty factor in the Eikonal equation.
We omit it here for the sake of simplicity.}.
In \cite{Hughes_2002} the equations for ``the potential'' $\varphi$ were derived based
on a physical analogy:
{\it ``Pedestrians have a common sense of the task (called potential) they face to reach
their common destination such that any two individuals at different locations having the same
potential would see no advantage to either of exchanging places.
There is no perceived advantage of moving along a line of constant potential. Thus, the motion
of any pedestrian is in the direction perpendicular to the potential [level sets].''}  
However, this suggestive analogy is somewhat misleading. 
After all, the gradient of the potential indicates the direction of forces/accelerations, while we are choosing the direction for the velocity.

Hughes further claims that
{\it ``Many physical situations are potentially both anisotropic and inhomogeneous. Fortunately,
[this] does not change the derivation
of the equations of motion and 
[the equivalent of $\bV(\x,t) = - v^2(\rho(\x,t)) \nabla \varphi(\x,t)$]
with [the speed] redefined for the new arguments,
remains valid.''} 
But it is easy to see that this is not quite correct: in the general anisotropic setting, the optimal direction can be different from $\frac{-\nabla \varphi(\x,t)}{|\nabla \varphi(\x,t)|}$.
(Consider the example of someone trying to swim to the opposite shore of a river of finite width and infinite length. Then, by symmetry, $\nabla\varphi$ points directly across the river from one shore to the other. However, the optimal path of the swimmer is to move in a diagonal direction, with the flow of the river, rather than fighting against the flow to go straight across. See also a detailed discussion and Figure \ref{fig:SwimmerRiver} in section \ref{ss:SpeedProfiles_Inter}.) Even though anisotropic problems were not considered by Hughes in \cite{Hughes_2002}, we believe that this comment may have 
led to the misidentification of optimal directions in several subsequent papers.


\section{Optimal directions and velocity profiles.}
\label{s:SpeedProfiles}
To evolve the crowd density $\rho$, it might seem natural to require that
the velocity field $\bV$ should be uniquely defined for all $(\x,t).$
However, the $\argmax$ approach for choosing the optimal direction of motion assumes that $\nabla \varphi$ is well-defined.
Unfortunately, equation \eqref{eq:HJB_general} typically does not have a smooth/classical solution, and admits infinitely many weak solutions.
Additional criteria based on test functions were introduced by Crandall and Lions in \cite{crandall1983viscosity} to select the unique {\em viscosity solution}, coinciding with the value function of the corresponding control problem.   For starting points where  $\nabla \varphi$ is undefined,
it is easy to see that optimal trajectories are usually non-unique.  Suppose $\{\tilde{\x}_n\}$ and $\{\hat{\x}_n\}$ are two sequences converging to the same point $\x$,
but $\tilde{\bp} = \lim_{n \to \infty}  \nabla \varphi(\tilde{\x}_n)$
is different from
$\hat{\bp} = \lim_{n \to \infty}  \nabla \varphi(\hat{\x}_n).$
Then it is clear that the optimal direction of motion starting from $\x$ will depend on whether we use $\tilde{\bp}$ or $\hat{\bp}$ in place of
$\nabla \varphi$ in equation \eqref{eq:HJB_general}.  Luckily, mild technical conditions on $\bv$
guarantee that the viscosity solution is Lipschitz-continuous (and, as a result, differentiable almost everywhere).  Since the crowd density $\rho$ is assumed to be in $L^1$,
it is sufficient to have $\bV$ uniquely defined almost everywhere, making the above complications less relevant.

Wherever $\varphi$ is differentiable,
the optimization in \eqref{eq:HJB_general} has a simple geometric structure: 
if we view $\bn = \frac{\nabla \varphi(\x)}{|\nabla \varphi(\x)|}$ as a known/fixed unit vector,
$\bu \in \U$ is simply chosen to maximize the projection of $\bv(\x, \bu)$ onto $(-\bn)$.
Two important consequences of this are related to the shape of the {\em velocity profile}
$\V(\x) = \{ \bv(\x,\bu) \, | \, \bu \in \U \}$:\\
1. The argmax in \eqref{eq:HJB_general} is a singleton if and only if the curve $\V(\x)$ is strictly convex.\\
2. If the curve $\V(\x)$ is smooth, each $\bu \in \U$ is optimal for at most one unit vector $\bn$.\\
(The latter property becomes ``if and only if'' for convex velocity profiles.)

In the isotropic case $\bv(\x, \bu) = v(\x) \bu$ and the velocity profile $\V(\x)$ is just a circle centered at the origin.  Thus, the optimal direction of motion is just $\bu_*(\x) = -\bn = \frac{-\nabla \varphi(\x)}{|\nabla \varphi(\x)|}$, as noted in the previous section for the Eikonal equation.
But as shown in Figure \ref{fig:InterCrowdSpeedProfiles}, this is not the case for general velocity profiles.

The lack of strict convexity for $\V(\x)$ implies that, for some $\bn^{\text{bad}}$ there are at least two different optimal directions of motion; e.g., see Figure \ref{fig:sub4}.
Unlike the  $\hat{\bp} \neq \tilde{\bp}$ scenario discussed earlier, this presents a significant problem since this $\bn^{\text{bad}}$ might be in effect on a large part of  $\domain$.

In subsections \ref{ss:SpeedProfiles_Inter}-\ref{ss:SpeedProfiles_Intra} we focus on two sources of anisotropy in pedestrian path planning.
The first is to consider the case of two pedestrian crowds ($A$ and $B$) moving towards different exits, with their speed depending on the angle between the directions of motion of the two crowds ($\bu_A$ and $\bu_B$). Examples of such models can be found in \cite{Jiang_2009, xiong2011, xiong2011lanes, jiang2012numerical, jiang2015higher, hanseler2015}.
Another source of anisotropy is to have the velocity of the pedestrians depend on the density of the crowd in a non-local, non-radially symmetric way, as in \cite{Cristiani_2015}.


\subsection{Inter-crowd anisotropy.}
\label{ss:SpeedProfiles_Inter}
Hughes' paper used isotropic pedestrian speed functions 
based on
the classical Lighthill and Whitham  model for vehicular flow\cite{Lighthill317}; e.g.,
$v(\rho,\bu) = \bar{v} \left(1 - \rho/\rho_{\text{max}}\right),$
where $\bar{v}$ and $\rho_{\text{max}}$ are positive constants.
  Other monotone non-increasing functions of $\rho$ have been similarly explored; e.g.,
the control-theoretic reinterpretation of Hughes' work in \cite{Huang_2009} used
a version of Drake's fundamental diagram \cite{drake1967}
\begin{equation}
\label{Huang_speed}
v(\rho,\bu) = \bar{v} e^{-\alpha\rho^2}.
\end{equation}

Hughes also considered multi-crowd scenarios, with each crowd distinguished by its desired destination and a different value of $\bar{v}$.
However, in his models the density induced slow-downs were due to the {\em total} density only.  E.g., given two crowds $A$ and $B$,
their speeds are specified by
$v^c(\rho,\bu) = \bar{v}^c \left(1 - \rho/\rho_{\text{max}}\right),$ where $c \in \{A,B\}$ and $\rho = \rhoA + \rhoB$.
This yields a total of four PDEs for $(\rhoA, \varphiA, \rhoB, \varphiB)$, with the coupling between the crowds based on $\rho$ only.
Importantly, these models are still isotropic and the equations for both $\varphiA$ and $\varphiB$ are Eikonal;
i.e., the differences in the directions chosen by two crowds have no effect on their 
speeds even when they are passing through  each other.

More recent work by Jiang et al. \cite{Jiang_2009} aimed to penalize for such directional differences explicitly by extending Formula \eqref{Huang_speed}.
The same approach was later adopted in \cite{xiong2011, xiong2011lanes, jiang2012numerical, jiang2015higher} and served as a background in \cite{jiang2017navigation_field}.
Assuming that $\buA$ and $\buB$ are the crowds' respective directions and the angle between them is $\psi = \arccos \left(\buA \cdot \buB \right)$, the ``disagreement penalty'' factor is defined as
\begin{equation}
\label{eq:square_penalty}
f(\bar{\rho}, \psi)  = e^{-\beta \left(1-\cos \psi \right)\bar{\rho}^2},
\end{equation}
where $\bar{\rho}$ is the density of the ``other'' crowd.
With this notation the new crowd speeds are defined as
\begin{equation}
\label{eq:penalized_speeds}
\vA \left(\rhoA,\rhoB,\buA,\buB \right)  = \bar{v} e^{-\alpha \left(\rhoA+\rhoB \right)^2} f(\rhoB, \psi)
\; \text{ and }  \;
\vB \left(\rhoA,\rhoB,\buA,\buB \right)  = \bar{v} e^{-\alpha \left(\rhoA+\rhoB \right)^2} f(\rhoA, \psi)
\end{equation}
for motion in the respective crowd directions $\buA$ and $\buB$; see the ``geometric velocity'' definition \eqref{eq:geom_velocity}.
The idea is that $f=1$ when both crowds move in unison, and $f$ monotonically decreases with both $\bar{\rho}$ and $\psi$, reaching the minimum when crowds move in opposite directions.
The values of parameters $\bar{v} = 1.034$, $\alpha = 0.075$, and $\beta = 0.019$  were
selected based on the physical experiments described in \cite{Wong_2010} (crowds of pedestrians asked to pass through each other at oblique angles with their aggregate speeds determined by video cameras).

\begin{figure}
\centering
\begin{subfigure}[t]{.4\textwidth}
  \centering
  \includegraphics[width=\linewidth]{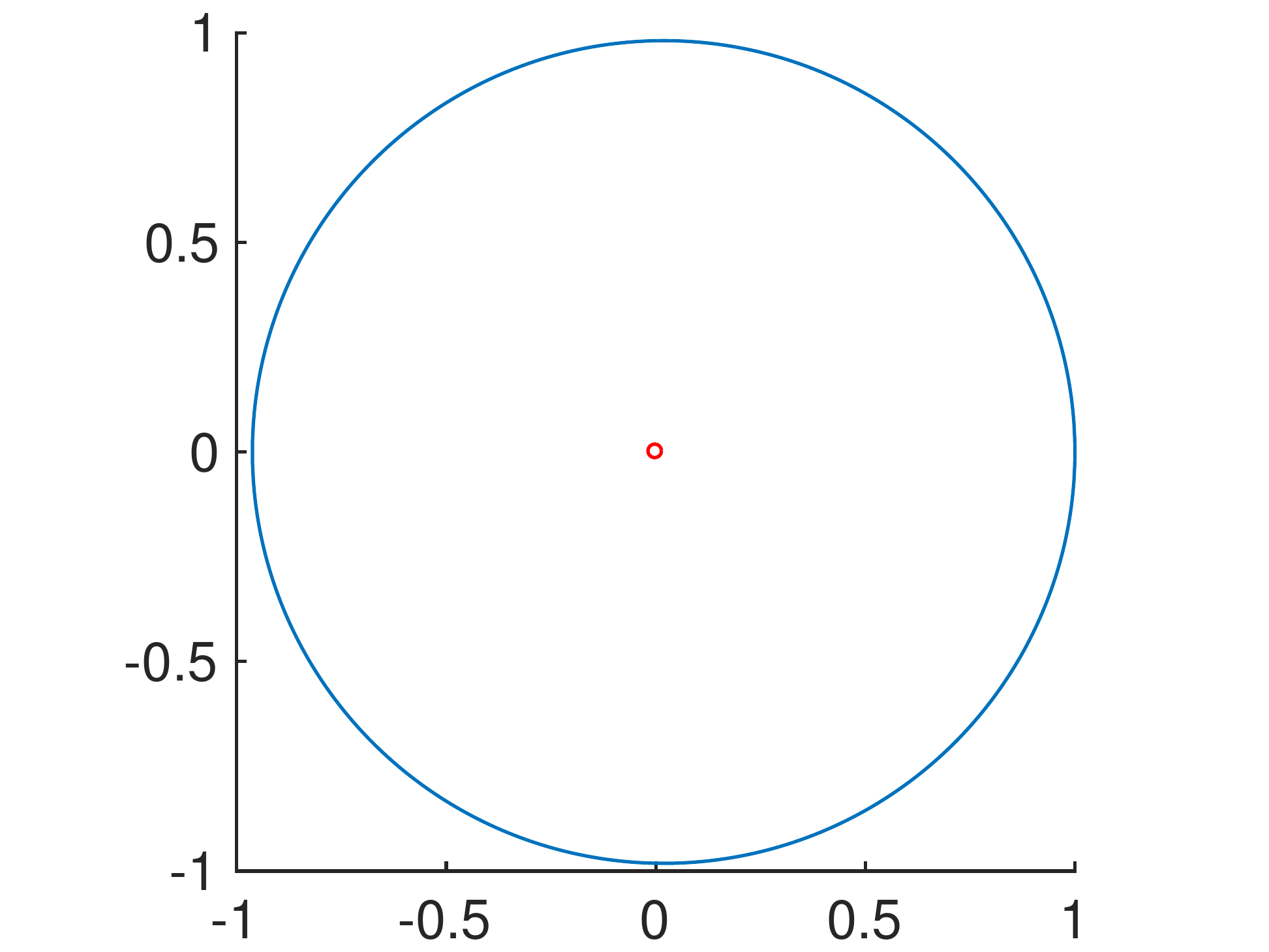}
  \caption{$\beta = 0.019$, $\rhobar = 1$ \newline Convex velocity profile}
  \label{fig:sub1}
\end{subfigure}
\begin{subfigure}[t]{.4\textwidth}
  \centering
  \includegraphics[width=\linewidth]{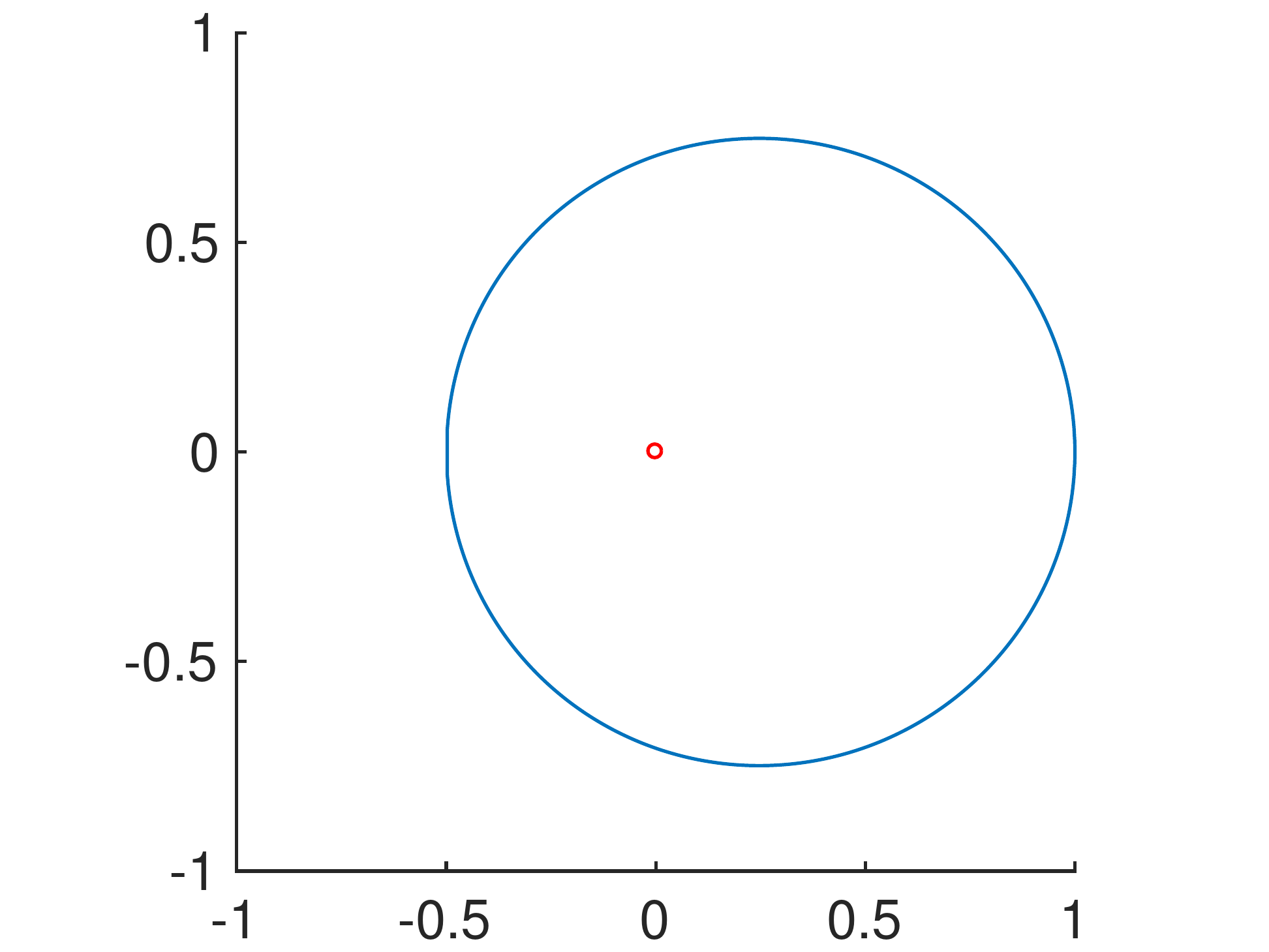}
  \caption{$\beta = 0.347$, $\rhobar = 1$ \newline Convex velocity profile}
  \label{fig:sub2}
\end{subfigure}
\begin{subfigure}[t]{.4\textwidth}
  \centering
  \includegraphics[width=\linewidth]{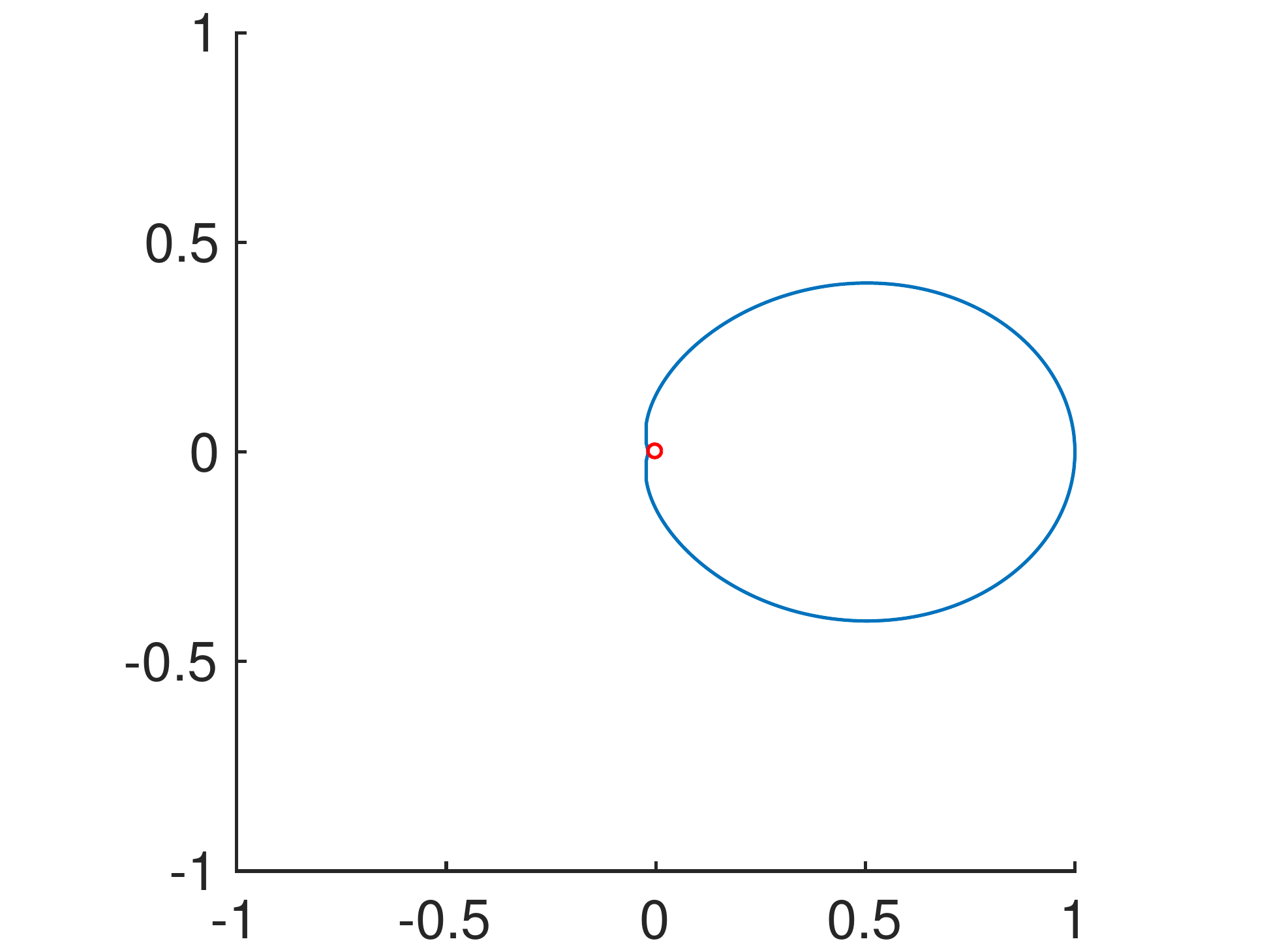}
  \caption{$\beta = 0.347$, $\rhobar = 2.4$ \newline Non-convex velocity profile}
  \label{fig:sub3}
\end{subfigure}
\begin{subfigure}[t]{.4\textwidth}
  \centering
  \includegraphics[width=\linewidth]{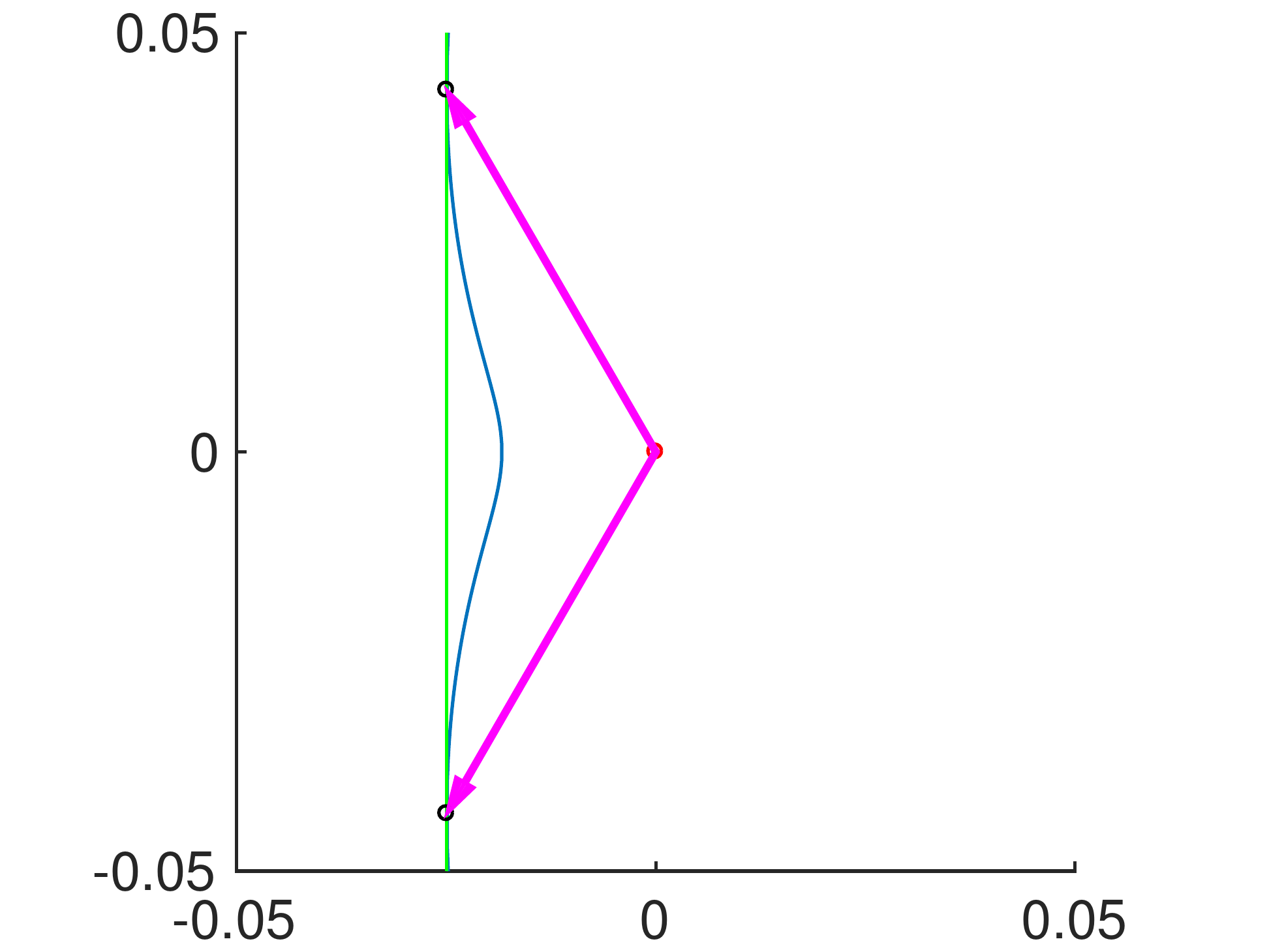}
  \caption{Zoomed version of $\beta = 0.347$, $\rhobar = 2.4$ to show non-convexity}
  \label{fig:sub4}
\end{subfigure}
\caption{Anisotropic velocity profiles for Crowd A at various parameter values.  We assume $\buB_* = (1,0).$}
\label{fig:InterCrowdSpeedProfiles}
\end{figure}

The dependence on $\psi$ clearly introduces a game-theoretic aspect: both crowds choose their directions simultaneously and their decisions affect each other immediately. Significant technical complications resulting from this strong coupling are discussed in section \ref{ss:Nash_Inter}.  But for now we note that even if $\buB_*$ is already chosen and known, the task of selecting an optimal $\buA_*$ is inherently anisotropic: the value function $\varphiA$ satisfies a version of the Hamilton-Jacobi equation \eqref{eq:HJB_general}
$$
\max\limits_{\buA \in \U}
\left\{-\nabla \varphiA(\x,t) \cdot \bvA\left(\rhoA(\x,t), \rhoB(\x,t), \buA, \buB_* \right) \right\} \; = \; 1.
$$
Figure \ref{fig:InterCrowdSpeedProfiles} illustrates the anisotropy in velocity profiles of $\varphiA$.
In \cite{Jiang_2009} and the more recent extensions \cite{xiong2011, xiong2011lanes, jiang2012numerical, jiang2015higher, jiang2017navigation_field} this anisotropic effect was ignored 
partly due to rewriting the above equation in ``Eikonal-type'' notation and partly due to the interpretation  of $\varphiA$ as a ``potential'', mirroring the original model by Hughes.  This led them to take $\frac{-\nabla \varphiA}{|\nabla \varphiA|}$ as the ``optimal'' direction, which in general can be quite different from the correct $\buA_*$.

However, it is worth noting that at these parameter values, the problem is very nearly isotropic at lower densities. At $\beta = 0.019$ and $\rhoB = 1$, we found the maximum angular difference between $\frac{-\nabla\varphi}{|\nabla\varphi |}$ and $\buA_*$ to be less than 0.02 radians
(See Figure \ref{fig:sub1} for a plot of the nearly circular velocity profile).
This might explain why this discrepancy was not noticed in the numerical experiments reported in \cite{Jiang_2009, xiong2011, xiong2011lanes, jiang2012numerical, jiang2015higher, jiang2017navigation_field}.

\begin{figure}
\centering
\begin{subfigure}[t]{.4\textwidth}
  \centering
  \includegraphics[width=\linewidth]{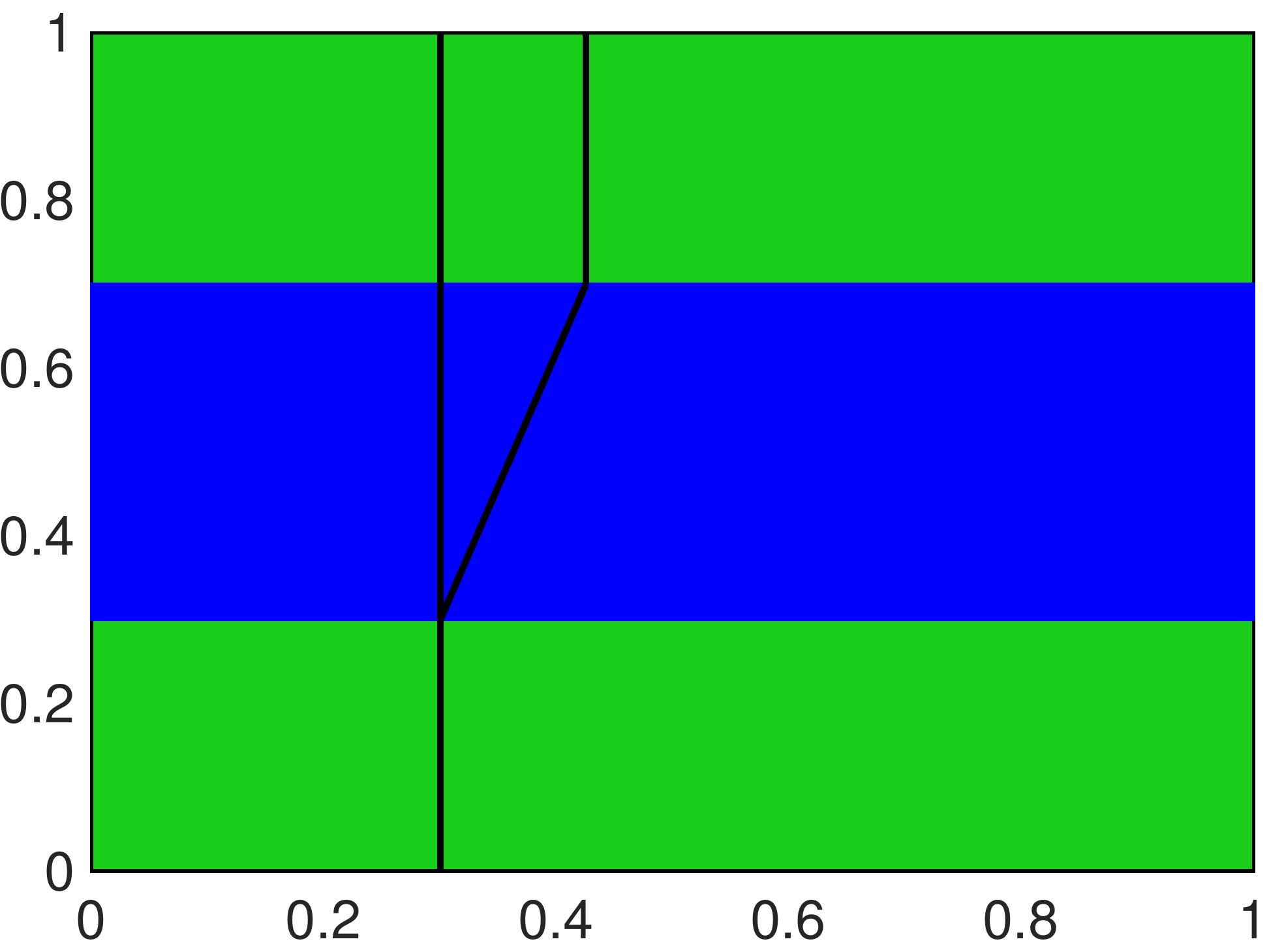}
  \caption{}
  \label{fig:SwimmerRiverSub}
\end{subfigure}
\begin{subfigure}[t]{.4\textwidth}
  \centering
  \includegraphics[width=\linewidth]{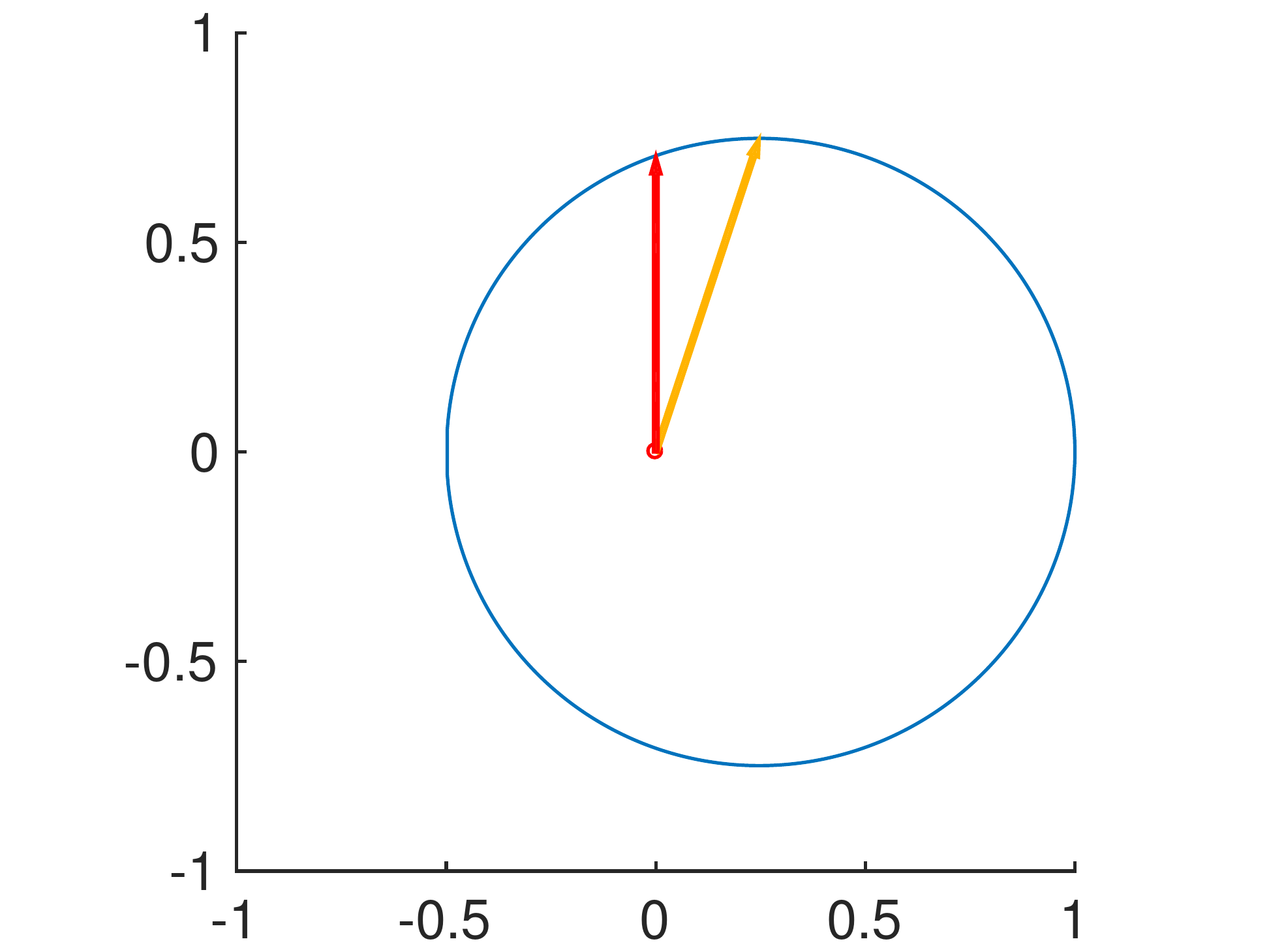}
  \caption{}
  \label{fig:SwimmerSpeedProfile}
\end{subfigure}
\caption{Optimal versus gradient path planning in subfigure (a). Corresponding velocity profile while ``in the river'' in subfigure (b). }
\label{fig:SwimmerRiver}
\end{figure}

To highlight how the optimal direction is in fact $\buA_*$ and not $\frac{-\nabla \varphiA}{|\nabla \varphiA|}$, let us consider an example with
Crowd $A$ consisting of a single pedestrian (i.e., $\rhoA = 0$ everywhere).
Let Crowd $B$ form a ``river'' of infinite length moving to the right, with density $\rhoB = 1$ for $0.3 \le y \le 0.7$, and $\rhoB = 0$ everywhere else.
Then let us consider the path planning of our single pedestrian of Crowd $A$ starting at $x=0.3$ and $y=0$, trying to reach the upper boundary $y=1$ as a target.
Then $\frac{-\nabla \varphiA}{|\nabla \varphiA|}$ always points directly upward.
However, the optimal path is actually to move diagonally with crowd $B$ while in the ``river''.
For parameter values of $\alpha = 0.075$ and $\beta = 0.347$, the exit time for the gradient path planning is $T = 1.21$.
But with the correct $\buA_*$ path planning, the exit time is instead $T = 1.1762$; see Figure \ref{fig:SwimmerRiver}. For comparison, the individual's velocity profile while in the ``river'' can be found in Figure \ref{fig:SwimmerSpeedProfile}. Here we can see that the largest vertical component of the speed corresponds to moving diagonally in the orange direction, the direction of the individual's path through the ``river'' in Figure \ref{fig:SwimmerRiverSub}.
In particular, it has a larger vertical component than would result from moving strictly up (see the red arrow in Figure \ref{fig:SwimmerSpeedProfile}).

We note that the parameter values used in \cite{Jiang_2009} and subsequent publications seem highly unrealistic despite the experimental data used for calibration in \cite{Wong_2010}.
For instance, consider two crowds with densities $\rhoA = \rhoB = 1$ person per square meter.
For convenience, we will define $\xi = f(1,0) - f(1,\pi) = 1-e^{-2\beta}$, which gives the
speed reduction resulting from these crowds moving in opposite directions rather than in unison.
At $\beta = 0.019$, if these crowds move in opposite directions, they are only slowed down  by $\xi \approx 4.7\%$ of the full speed achievable when they move in unison.
To provide some perspective, the classical estimate by Navin and Wheeler \cite{NavinWheeler} used to calibrate many traffic engineering models
yields\footnote{
Navin and Wheeler \cite{NavinWheeler} have examined the case of two crowds heading in opposite directions with $\rhoA + \rhoB = 1$ and different ratios of $\nu = \rhoA/(\rhoA+\rhoB)$.
The case of $\nu  = 1$ is equivalent to 2 crowds moving in unison, while $\nu  = 1/2$ yields $\Psi = \pi$ and $\rhoA = \rhoB = 1/2.$
In the latter case they have observed that the speeds were reduced by $4\%$.}
a significantly higher $\beta = 0. 078.$
For our example, with $\rhoA = \rhoB =  1$, this corresponds to a slow down of $\xi \approx 14.5\%$ relative to moving in unison.
Even this last estimate appears rather conservative: anyone who has ever tried to catch a subway train during rush hour would
attest that a speed reduction of $\xi \in [30\%,50\%]$ (corresponding to $\beta \in [0.178,0.347]$) would be much closer to their experience.

For high enough values of $\rhoB$ or $\beta$, the velocity profile of crowd A can actually become non-convex.
This scenario is highlighted in Figure \ref{fig:sub4}. If a pedestrian tries to maximize the component of $\buA$ in the $-x$ direction, then both choices of $\buA$ denoted by the magenta arrows are optimal.

When looking at the velocity profile $\V$,
we are concerned with only the smoothness and convexity of the profile, and not the scaling.
Thus we need only look at the ``disagreement penalty'' factor $f(\rhoB, \psi)$.  Assuming that $\rhoB$ is fixed, we will further refer to it as $f(\psi)$ to simplify the notation.
For a fixed $x$, we only need to check that the polar curve $\textbf{f}(\psi) = (f(\psi)\cos(\psi),f(\psi)\sin(\psi))$ is strictly convex, where $\psi$ is the angle between $\buA$ and $\buB$.
If we consider its tangent vector $\textbf{T}(\psi) = \dot{\textbf{f}}(\psi)$,
we can then write the derivative of the tangent vector as $\textbf{T}'(\psi) = k(\psi)\textbf{N}(\psi)$, where $k$ is the curvature, and $\textbf{N}$ is the unit normal vector.
The sufficient condition for strict convexity is that $\textbf{T}(\psi) \times \textbf{N}(\psi)$ does not change sign.
This reduces to a simple condition on $f$ and its derivatives:
\begin{equation}\label{eq:intercrowd_convexity}
f(\psi)^2 + 2f'(\psi)^2 - f(\psi)f''(\psi) > 0, \quad \forall \psi
\end{equation}

For the particular model \eqref{eq:square_penalty}, this condition simplifies to:
$$ 1+\beta\rhobar^2 \cos(\psi) + \beta\rhobar^2\sin^2(\psi) > 0. $$

We see that the left hand side is minimized when $\psi = \pi$.
The sufficient condition for strict convexity then becomes
$$ \beta\rhobar^2 < 1. $$

At the parameter values used in \cite{Jiang_2009}, we see that we only lose strict convexity at densities exceeding $\rhobar \approx 7.25$ people per square meter -- a threshold which is too high to be relevant in practice. However, as discussed previously, their value of $\beta$ is rather small and with $\beta = 0.078$ (consistent with Navin and Wheeler's data \cite{NavinWheeler}), we instead get that we can lose strict convexity when the density of either crowd exceeds a more physically relevant threshold of $\rhobar = \sqrt{1/0.078} \approx 3.57.$

We summarize the results of using different values of $\beta$ in Table \ref{table:RhoTable}.
The first two columns are for $\beta$ values consistent with Wong et al. \cite{Wong_2010} and Navin et al. \cite{NavinWheeler} respectively.
The third and fourth columns correspond to a slowdown of $\xi = 0.3$ and $0.5$ respectively.
In each case we list the critical density $\hat{\rho}$ at which the velocity profile loses strict convexity.
We note that different ``disagreement penalty factors'' can be studied similarly.
For example, a model considered in \cite{hanseler2015} is based on using $\rhobar$ rather than $\rhobar^2$; i.e.,
\begin{equation}
\label{eq:non_square_penalty}
f(\bar{\rho}, \psi)  = e^{-\beta \left(1-\cos \psi \right)\rhobar},
\end{equation}
with the advantage that it scales more consistently for pairwise interactions of multiple crowds.
The critical $\hat{\rho}$ values for this alternative model are reported in the last row of Table \ref{table:RhoTable}.

\begin{table}
\caption{Critical densities for two disagreement penalty models and various $\beta$ values.}
\label{table:RhoTable}
\centerline{
\begin{tabular}{|c|c|c|c|c|}
\hline
 $\beta$ & 0.019 & 0.078 & 0.178 & 0.347 \\
\hline
$\xi$ & 0.037 & 0.144 & 0.3 & 0.5 \\
\hline\hline
critical density $\hat{\rho}$ for the $\beta\rhobar^2$ model & 7.25 & 3.58 & 2.37 & 1.70 \\
\hline\hline
critical density $\hat{\rho}$ for the $\beta\rhobar$ model  & 52.6 & 12.8 & 5.62 & 2.88 \\
\hline\end{tabular}
}
\end{table}

\subsection{Intra-crowd anisotropy.}
\label{ss:SpeedProfiles_Intra}
Up until now we have only discussed models with localized interactions of pedestrians.
But it is often reasonable to let the velocity $\bV(\x,t)$ depend on values of $\rho(\cdot, t)$ in some region containing $\x$, rather than on a local $\rho(\x,t)$ only.
For multiple crowds, such non-local models present additional analytic challenges even if the direction field is specified in advance instead of chosen by individual pedestrians \cite{colombo2012nonlocal}.  If individual choices are also included, non-local interactions give rise to anisotropy. In this section we show that this anisotropy may result in model inconsistency even for a single crowd.

\begin{figure}
\centerline{
\begin{tikzpicture}[scale=3]
	\begin{scope}[thick]
	\draw (0,0) -- (0.996,0.0872);
	\draw (0,0) -- (-0.996,0.0872);
	\draw[->,  very thick] (0,0) -- (0,0.7);
	\draw (0.996,0.0872) arc (5:175:1);
	\end{scope}[semithick]
	\draw[thin] ( 0.0996, 0.00872) arc (5:175:0.1);
	\draw (0.1,0.15) node {$\alpha$};
	\draw (0.7,0) node {$R$};
	\draw (-0.05,-0.05) node {$\x$};
	\draw (0.1,0.7) node {$\bu$};
	\draw (-0.5,0.6) node {$S(\x,\bu)$};
\end{tikzpicture}
}
\caption{An example of sensory region $S(\x,\bu)$ for intra-crowd anisotropy.}
\label{fig:SensoryRegion}
\end{figure}
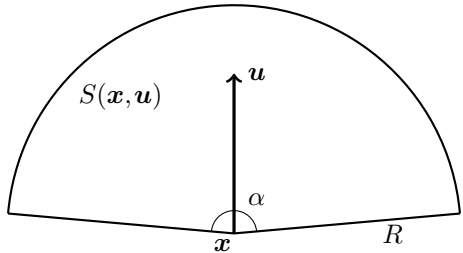

In the paper by Cristiani, Priuli, and Tosin \cite{Cristiani_2015},
the velocity is no longer chosen directly. Instead, it is the sum of a chosen behavioral velocity $\bu \in S^1$, and an interaction velocity perturbation $\w_i$ that is determined based on the density of pedestrians inside some sensory region $S(\x,\bu)$. This sensory region is usually modeled by a sector (of angle $\alpha$ and radius $R$) centered at $\x$ and with the orientation determined by $\bu$; see Figure \ref{fig:SensoryRegion}.
The key idea is that each pedestrian is only influenced by others roughly ahead of her, and not the ones behind. The model developed in \cite{Cristiani_2015} defines the interaction velocity $\w_i$ as follows:
\begin{equation}\label{eq:IntraCrowdSpeed}
\w_i
= \int_{S(\x,\bu) \cap \Omega} \mathcal{F}(\y-\x)\rho(t,\x) \, d\y;
\qquad \qquad
\mathcal{F}(\br) =
\begin{cases}
\frac{-F\br}{|\br|^2}, &|\br| > \frac{1}{C}, \\
\frac{-FC \br}{|\br|}, &|\br| \le \frac{1}{C},
\end{cases}
\end{equation}
where $F$ and $C$ are positive constants, with the latter functioning as a ``cutoff''.
Since the interaction velocity perturbation has a non-local dependence on the current $\rho$, it would be accurate to denote it
$\w_i \left[ \rho(\cdot, t) \right] (\x, \bu)$, but we will at times abbreviate this as $\w_i(\x,\bu)$ or even $\w_i$. 
The overall velocity is 
\begin{equation}\label{eq:nonLocal_velocity}
\bv \left[ \rho(\cdot, t) \right](\x,\bu) \; = \; \bu \, + \, \w_i\left[ \rho(\cdot, t) \right](\x,\bu),
\end{equation} 
which is a version of formula \eqref{eq:additive_velocity}.
Modulo this change from $\bv \left( \rho(\x, t), \bu \right)$ to $\bv \left[ \rho(\cdot, t) \right](\x,\bu)$, the functions $\rho$ and $\varphi$ satisfy the same system of PDEs \eqref{eq:strobo_full}.

We will now investigate the convexity of the velocity profiles for the case of intra-crowd anisotropy.
Since this model depends non-locally on the density, we will only analyze this under some assumptions on the density function.

\subsubsection{Linear Density.}
First, let us consider the case of a density distribution that is linear in $\x$ at some fixed time $t$.
WLOG we may assume that this density distribution is of the form:
$$ \rho(x,y) = \rho_0 + \rho_x x $$
In this case, we can evaluate the integral found in \eqref{eq:IntraCrowdSpeed}, which is done in Appendix \ref{ss:Appendix1}.

\begin{figure}
\centerline{
\begin{subfigure}[t]{0.4\textwidth}
\includegraphics[width=\linewidth]{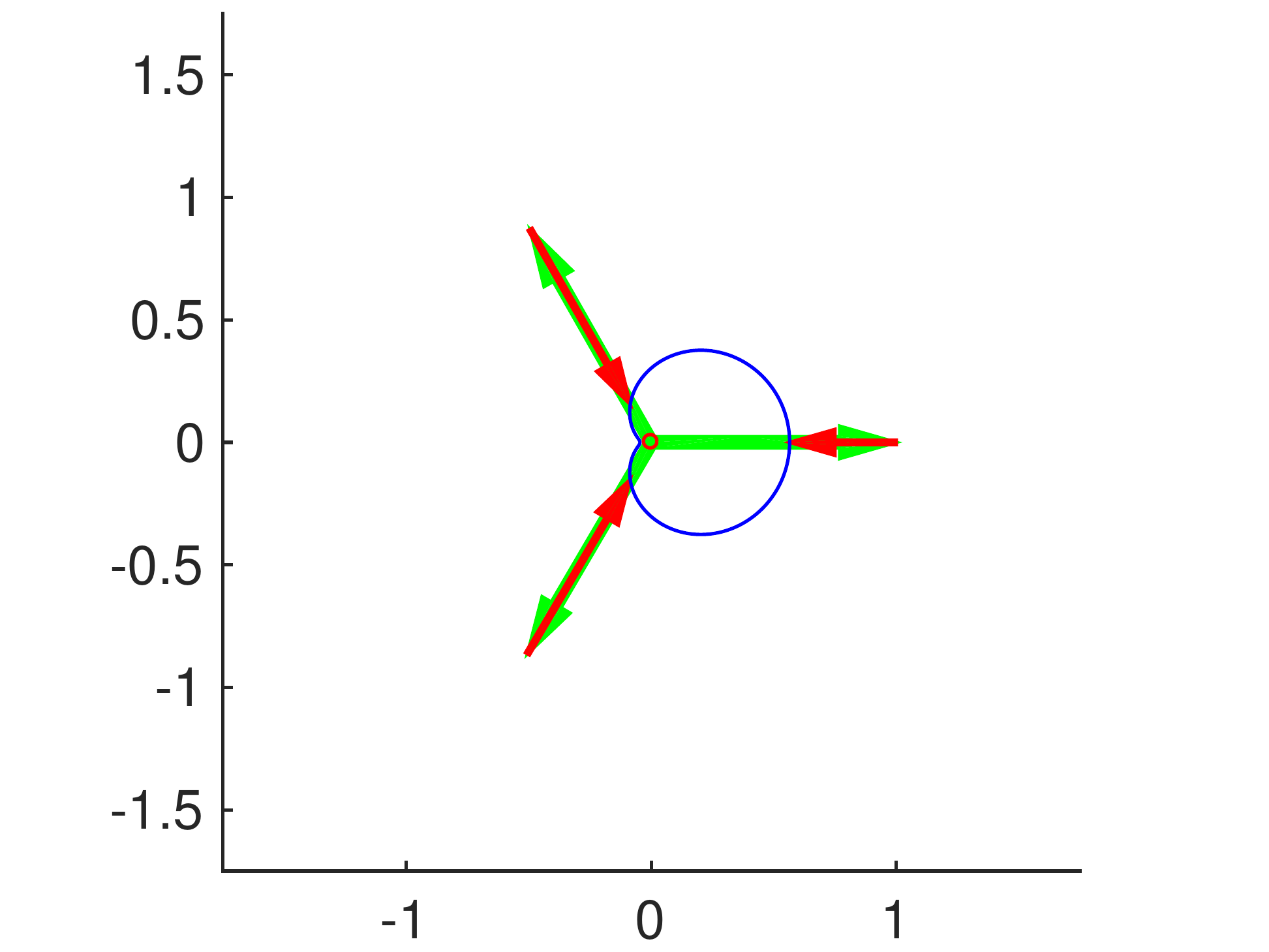}
\caption{}
\label{fig:IntraCrowdSpeedProfile}
\end{subfigure}
\begin{subfigure}[t]{0.4\textwidth}
\includegraphics[width=\linewidth]{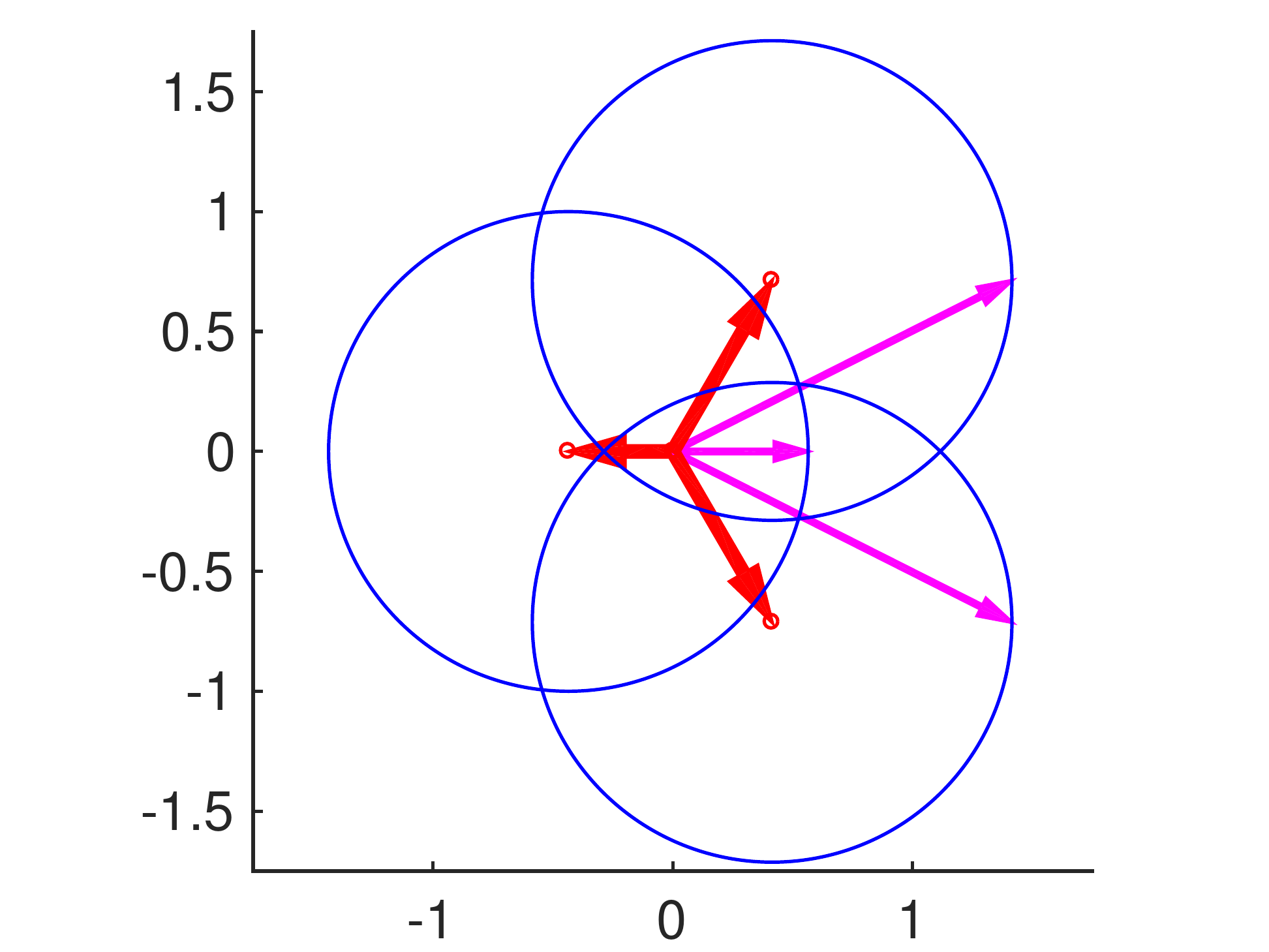}
\caption{}
\label{fig:ShiftedCircle}
\end{subfigure}
}
\caption{Left: Non-convex velocity profile for intra-crowd anisotropy and linear density with $\alpha = 20^\circ, F=1, R=1, \rho_0 = 2, \rho_x = -3/2$ on the left.  Right: three ``shifted circle'' versions corresponding to different/fixed $\bw_i$ choices.  Magenta arrows show the resulting velocity corresponding to $\bu = (1,0)$ in each profile.  One of them coincides with the correct resulting velocity in subfigure (a).}
\end{figure}

We note that for very high densities, the velocity profile may not actually contain the origin, as the interaction velocity can dominate the behavioral velocity.
We show in Appendix \ref{ss:Appendix2} that as long as the density is linear in $\x$ and the velocity profile contains the origin, then we have convexity for all $\alpha$ such that $\alpha > 3\sin(\alpha)$.
This condition is satisfied for all $\alpha > \approx 130.57^\circ$, including the value of $\alpha = 170^\circ$ used by Cristiani et al \cite{Cristiani_2015}.
Otherwise, it may be possible to find a linear density function, and values for the parameters such that the velocity profile will not remain convex.
For example, consider the case of $\alpha = 20^\circ, F=1, R=1$, $\rho_0 = 2$, and $\rho_x = -3/2$. A plot of this velocity profile can be found in Figure \ref{fig:IntraCrowdSpeedProfile}.

This linear density example is also good for illustrating another (general) pitfall in anisotropic modeling.
In \cite{Cristiani_2015} the dynamics are given in the form:
$$
\y'(t) = \bu(t) + \bw_i[\rho(\tau,\cdot )](t)
$$
However, this misleadingly suggests that $\bw_i$ is independent of the choice of $\bu$.
To illustrate this point, we show the correctly computed velocity profile in Figure \ref{fig:IntraCrowdSpeedProfile} with
three different versions of $\bu$  in green and the corresponding  $\bw_i$ vectors in red.
Choosing any one of  these $\bw_i$'s to be used {\em for all} $\bu$ directions will yield a ``shifted circle'' velocity profile with very different resulting velocities; see Figure \ref{fig:ShiftedCircle}.
Moreover, it is not clear which of these shifts would be a reasonable choice.
This limits the validity of the simplified equation proposed in  \cite{Cristiani_2015}
$$
|\nabla \varphi| \, - \, \nabla \varphi \cdot  \bw_i[\rho(\tau,\cdot )] \; = \; 1
$$
to radially symmetric sensing regions (i.e., with $\alpha = 360^\circ$ only).
Based on personal communications with the authors, the actual numerical tests in \cite{Cristiani_2015} have used a correct semi-Lagrangian discretization of the general HJB equation (\ref{eq:strobo_full}b).


\subsubsection{Piecewise Constant Density.}
We also considered piecewise constant density functions of the form:
$$
\rho(x,y) =
\begin{cases}
\rho_1, &x > 0 \\
\rho_2, &x \le 0
\end{cases}
$$

\begin{figure}
\begin{subfigure}[t]{0.4\textwidth}
\includegraphics[width=\linewidth]{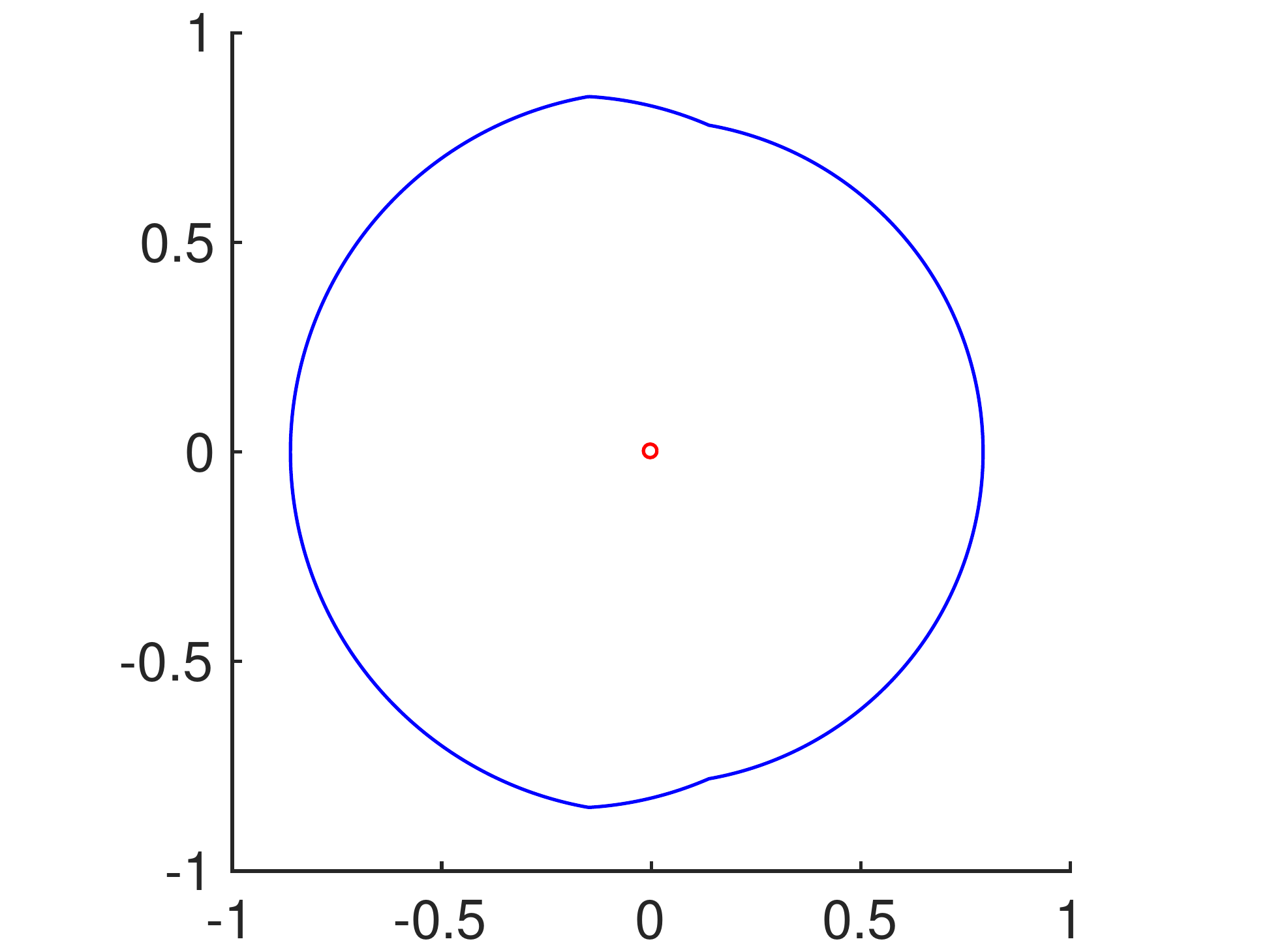}
\caption{}
\label{fig:discont1}
\end{subfigure}
\begin{subfigure}[t]{0.4\textwidth}
\includegraphics[width=\linewidth]{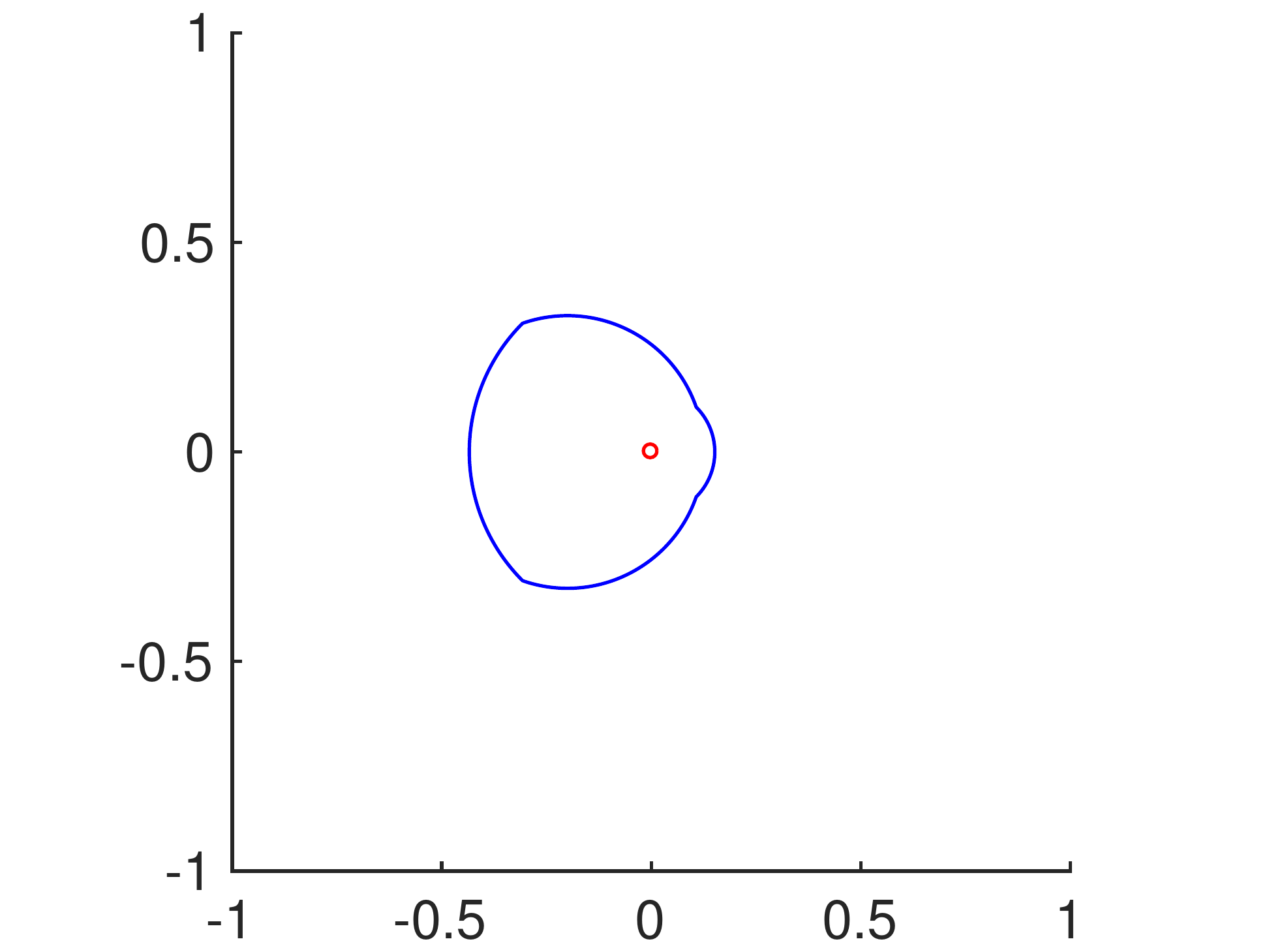}
\caption{}
\label{fig:discont2}
\end{subfigure}
\caption{Non-convex velocity profiles for discontinuous density: $F = 0.4$, $R = 1$, $\rho_1 = 1.5$, $\rho_2 = 1$, $\alpha = 20^\circ$ in (A) and the same except for $\alpha = 90^\circ$ in (B).}
\label{fig:discont}
\end{figure}

Here we find that both convexity and smoothness can also fail at $x=0$ for some parameter values, see the examples in Figure \ref{fig:discont}.

\section{Nash equilibria.}
\label{s:Nash}
Returning to the problem of two-crowd interactions with $\rho = (\rhoA,\rhoB),$ let us consider a speed function of the form:

\begin{equation*}
\vA (\rho,\buA) = f(\rhoB,\psi)
\end{equation*}
with the corresponding velocity
\begin{equation*}
\bvA (\x,\buA) = \buA f(\rhoB(\x),\psi)
\end{equation*}

where $\buA$ and $\buB$ are the directions of motion of each crowd and $\psi$ is the angle between them.

In this case, the decision-making of the two crowds is modeled as a non-zero-sum game.
We then look for a solution to the following coupled system of Hamilton-Jacobi-Isaacs (HJI) equations:
\begin{equation}
\max_{\buA \in S^1} \Bigl\{ (\nabla \varphiA(\x) \cdot (-\buA) \vA(\rhoA(\x),\buA,\buB)\Bigr\} = 1
\label{eq:HJI_1}
\end{equation}
\begin{equation}
\max_{\buB \in S^1} \Bigl\{ (\nabla \varphiB(\x) \cdot (-\buB) \vB(\rhoB(\x),\buB,\buA)\Bigr\} = 1
\label{eq:HJI_2}
\end{equation}
This notation is slightly misleading: each of the two PDEs involves maximizing in only one variable ($\buA$ or $\buB$), with the other one participating as a parameter.
Both maximizations are assumed to be performed simultaneously, yielding the optimal pair of direction fields $(\buA_*, \buB_*).$
A different pair of PDEs could be obtained by assuming that one of the crowds has an advantage of selecting its direction second (once the other crowd's direction is already known).
In general, the solution might change depending on which crowd we select to ``go first''.  This feature is disconcerting from a modeling perspective since the notion of ``going first'' is harder to motivate when the replanning is done continuously in time (with changing densities $\rhoA$ and $\rhoB$).  So, it is desirable to guarantee that the order of decision making is irrelevant, making the system (\ref{eq:HJI_1}, \ref{eq:HJI_2}) meaningful as written above.  The sufficient condition for this is the uniqueness of the two-player {\em Nash equilibrium} (NE), a concept which we review below for general games.


We assume that player $A$ chooses strategies $a \in \mathcal{A}$ to maximize $\JA(a,b)$ and player $B$ chooses strategies $b \in \mathcal{B}$ to maximize $\JB(a,b)$.
In the non-zero-sum setting (when $\JA \neq -\JB$) the well-posedness of the game is non-trivial even when $\mathcal{A}$ and $\mathcal{B}$ are discrete sets.
NE is defined as a pair of strategies $(a_*,b_*) \in \mathcal{A} \times \mathcal{B}$ such that there is never an incentive to change unilaterally. More precisely, $(a_*,b_*)$ is said to be NE if both of the following hold:
$$ \JA(a_*,b_*) \ge \JA(a,b_*), \quad \forall a \in \mathcal{A}, $$
$$ \JB(a_*,b_*) \ge \JB(a_*,b), \quad \forall b \in \mathcal{B}. $$
If the game has a unique NE, then we can always assume that both players make their decisions simultaneously -- if both players are rational, there is no longer any advantage to ``go second.''

The existence of NE for infinite continuous games has been proved under some very general assumptions, at least in the class of mixed policies \cite{glicksberg1952}. For zero-sum-games, the possible non-uniqueness is not a major issue and requires no coordination between the players.
If both $(a_*,b_*)$ and $(a_\sharp,b_\sharp)$ are NE, then so are $(a_*,b_\sharp)$ and $(a_\sharp,b_*)$; moreover, $\JA$ will have the same value for each of these pairs (and $\JB = -\JA$).  So, who goes first and which NE they aim for still does not matter.
However, none of this holds in the case of non-zero-sum games. In section \ref{ss:ExamplesNash} we will show examples where the value of multiple NE is different, players disagree about their relative merits, and choosing simultaneously without coordination can easily lead to a non-Nash outcome.

In the case of zero-sum HJI games, the notions of ``upper value'' and ``lower value'' are used to capture what happens when one of the players goes first.   A separate ``Isaacs criterion'' on the Hamiltonian ensures that these two values are always equal, enabling simultaneous decision-making by the players.
Additional mild technical conditions guarantee that a NE exists among {\em pure} strategies, making it unnecessary to consider any relaxed/``chattering'' controls \cite{bardicapuzzodolcetta}.  A simple modification of this argument yields the existence of a pure-strategy NE in (\ref{eq:HJI_1}, \ref{eq:HJI_2}) provided both velocity profiles are strictly convex.  However, to claim that this non-zero-sum system has a well-defined value, we will also need to know that this NE is {\em unique}.

We will focus on a fixed $\x \in \Omega$, assuming that $\rho$'s and $\varphi$'s are known for both crowds and writing  $\textbf{p}= \nabla\varphi_\textsc{a}(\textbf{x}), \; \textbf{q}=\nabla\varphi_\textsc{b}(\textbf{x})$ to simplify the notation:
$$
\max_{\buA \in S^1} \Bigl\{ \bp \cdot (-\buA) \,  \vA\left(\rhoB(\x),\buA,\buB\right)\Bigr\} = 1; \hfill
\qquad \qquad
\max_{\buB \in S^1} \Bigl\{ \bq \cdot (-\buB)  \, \vB\left(\rhoA(\x),\buA,\buB\right)\Bigr\} = 1.
$$
Furthermore, to focus on the geometric aspects we will use $\Delta$ to denote the angle between  $\bp$ and $\bq$,
$\psi$ for the angle between $\buA$ and $\buB$, and $(a,b)$ for the angles that $\buA$ and $\buB$ make with $(-\bp)$ and $(-\bq)$ respectively.

From the definition \eqref{eq:penalized_speeds},
we note that the anisotropic speeds $\vA$ and  $\vB$ have a common isotropic factor which does not affect the choice of optimal direction.
Denote $\fA(\psi) = f(\rhoB, \psi)$, $\fB(\psi) = f(\rhoA, \psi)$, recall that $\psi = a-b + \Delta$, and our new task is
to show the uniqueness of the NE for all $(\bp, \bq)$ in the equivalent game
\begin{equation}
\max_{a \in [0,2\pi ]} \Bigl\{\cos(a) \fA(a-b + \Delta)\Bigr\} = 1,
\hfill
\qquad \qquad
\max_{b \in [0,2\pi ]} \Bigl\{\cos(b) \fB(a-b + \Delta)\Bigr\} = 1.
\label{eq:game}
\end{equation}

We will assume that the velocity profiles for both crowds are smooth and strictly convex.
The strict convexity ensures that the $\argmax_{a \in [0,2\pi ]}$ is unique for every $b$,
making the ``best reply'' map $\PhiA(b)$ well-defined. Since $\fA$ is smooth, the first-order optimality condition
\begin{equation}
0 \; = \; \frac{d}{da} \Bigl[\cos(a) \fA(a-b + \Delta)\Bigr] = -\sin(a) \fA(a-b+\Delta) + \cos(a)\frac{\partial}{\partial a} \bigl[\fA(a-b+\Delta)\bigr]
\label{eq:opt_cond}
\end{equation}
should hold for $a = \PhiA(b).$
If we denote the right hand side as $g(a,b)$, then
$$
\frac{\partial g}{\partial a} \, = \, -\cos(a)\fA(\psi) -2\sin(a)\dfadpsi+\cos(a)\dfaddpsi;
\qquad
\frac{\partial g}{\partial b} \, = \, -\sin(a)\dfadpsi+\cos(a) \dfaddpsi.
$$
Noting that \eqref{eq:opt_cond} implies $\sin(a) = \cos(a)\frac{\dfadpsi}{\fA(\psi)},$  we can now simplify the derivative
\begin{equation}
\frac{d\PhiA}{db} \, = \, \frac{-g_b}{g_a} \, = \,
 \frac{-\left(\dfadpsi\right)^2 + \fA(\psi)\dfaddpsi}{-\fA(\psi)^2-2\left(\dfadpsi\right)^2 + \fA(\psi)\dfaddpsi}.
 \label{eq:dphi}
\end{equation}
A similar derivation can be repeated for the other crowd's best reply $\PhiB(a)$.
If we show that the map $a \mapsto \PhiA(\PhiB(a))$ is a contraction mapping,
the existence and uniqueness of NE becomes a simple consequence of the Banach fixed-point theorem.

Thus, it is sufficient to show that the derivative is bounded by 1, i.e.
$$
\left| \frac{d}{da} \Bigl[ \PhiA(\PhiB(a)) \Bigr] \right| < 1, \qquad \forall a \in [0,2\pi]
$$
It is generally easier to check the slightly stronger condition following from the Chain Rule that
$$
\left| \frac{d\PhiA}{db} \right| \left| \frac{d\PhiB}{da} \right| < 1, \qquad \forall (a,b) \in [0,2\pi] \times [0,2\pi]
$$
where the derivatives of the best reply maps can be found in Equation \eqref{eq:dphi}.

\subsection{NE in intercrowd models.}
\label{ss:Nash_Inter}
In \cite{Jiang_2009}, they assume a speed function for crowd $A$ of the form:
\begin{equation*}
\vA(\psi) = v_0e^{-\alpha(\left(\rhoA\right)^2+\left(\rhoB\right)^2)}e^{-\beta\left(\rhoB\right)^2(1-\cos(\psi))}
\end{equation*}
where $\rho$ is the density of the other crowd. (Crowd $B$'s speed function is defined analogously.)
Since we are interested only in the shape (rather than scaling) of the velocity profile, we shall just use:
\begin{equation*}
\fA(\psi) = e^{-\beta\left(\rhoB\right)^2(1-\cos(\psi))}
\end{equation*}
as it is simpler to work with.
For ease of notation, we shall also drop the superscripts on $\rho.$
We calculate that
\begin{equation*}
\dfadpsi = -\beta\rho^2\sin(\psi)e^{-\beta\rho^2(1-\cos(\psi))},
\end{equation*}
\begin{equation*}
\dfaddpsi = \beta^2\rho^4\sin(\psi)^2e^{-\beta\rho^2(1-\cos(\psi))} -\beta\rho^2\cos(\psi)e^{-\beta\rho^2(1-\cos(\psi))}.
\end{equation*}
Plugging these into Equation \eqref{eq:dphi}, we obtain
\begin{equation}
\left| \frac{d\PhiA}{db} \right| \, = \, \left| \frac{\beta\rho^2\cos(\psi)}{\beta\rho^2\cos(\psi)+\beta^2\rho^4\sin^2(\psi)+1} \right|.
\end{equation}

We see that the denominator is a quadratic polynomial in $\cos(\psi)$. If we assume the necessary and sufficient condition for strict convexity that $\beta\rho^2 < 1$, then it follows that the absolute value of this polynomial reaches a minimum of $1-\beta\rho^2$ at $\psi = \pi$. Similarly, the numerator is maximized in absolute value at $\psi = \pi$, so we see that the worst case scenario is that of $\psi = \pi$, where the derivative is bounded by:
\begin{equation}
\left| \frac{d\PhiA}{db} \right| \, < \frac{\beta(\rhoB)^2}{1-\beta(\rhoB)^2}.
\end{equation}
The same process works with the crowds switched.
So $\PhiA(\PhiB(a))$ is a contraction mapping (and NE is unique) whenever
\begin{equation}
\label{eq:Nash_uniqueness}
\frac{\beta(\rhoB)^2}{1-\beta(\rhoB)^2} \frac{\beta(\rhoA)^2}{1-\beta(\rhoA)^2} < 1
\qquad \text{or, equivalently,} \qquad
\beta(\rhoA)^2 + \beta(\rhoB)^2 < 1.
\end{equation}

\begin{figure}[h]
$
\arraycolsep=30pt
\begin{array}{cc}
\begin{tikzpicture}[scale=2]
	\begin{scope}[semithick]
	\draw (0,0) -- (1.5,0);
	\draw (0,0) -- (0,1.5);
	\draw (1,0) -- (1,1);
	\draw (0,1) -- (1,1);
	\draw (1.5,0) -- (1.45,-0.05);
	\draw (1.5,0) -- (1.45,0.05);
	\draw (0,1.5) -- (0.05,1.45);
	\draw (0,1.5) -- (-0.05,1.45);
	\draw (1,0) arc (0:90:1);
	\end{scope}[semithick]
	\draw (0.5,0.5) node {$1$};
	\draw (0.85,0.85) node {$2$};
	\draw (1.25,1.25) node {$3$};
	\draw (1,-0.15) node {$\hat{\rho}$};
	\draw (-0.1,1) node {$\hat{\rho}$};
	\draw (1.6,0) node {$\rhoA$};
	\draw (0,1.6) node {$\rhoB$};
\end{tikzpicture} &
\begin{tikzpicture}[scale=2]
	\begin{scope}[semithick]
	\draw (0,0) -- (1.5,0);
	\draw (0,0) -- (0,1.5);
	\draw (1,0) -- (1,1);
	\draw (0,1) -- (1,1);
	\draw (1.5,0) -- (1.45,-0.05);
	\draw (1.5,0) -- (1.45,0.05);
	\draw (0,1.5) -- (0.05,1.45);
	\draw (0,1.5) -- (-0.05,1.45);
	\draw (1,0) -- (0,1);
	\end{scope}[semithick]
	\draw (0.25,0.25) node {$1$};
	\draw (0.75,0.75) node {$2$};
	\draw (1.25,1.25) node {$3$};
	\draw (1,-0.15) node {$\hat{\rho}$};
	\draw (-0.1,1) node {$\hat{\rho}$};
	\draw (1.6,0) node {$\rhoA$};
	\draw (0,1.6) node {$\rhoB$};
\end{tikzpicture}\\
(a) & (b)
\end{array}
$
\caption{Critical density regions for the $\beta(\rhobar)^2$ model (subfigure a) and for the $\beta\rhobar$ model (subfigure b).
Provably unique NE in Region 1; multiple NE observed in Region 2; non-convex velocity profiles in Region 3.}
\label{fig:DensityRegions}
\end{figure}
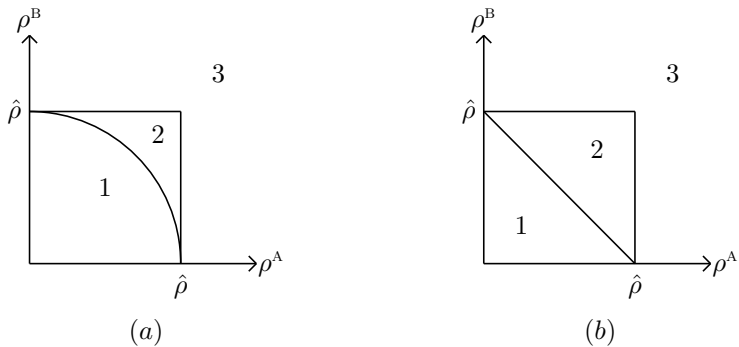

The models in \cite{Jiang_2009} assume a value of $\beta = 0.019$, which means that our condition for guaranteeing the uniqueness of NE becomes:
\begin{equation}
(\rhoA)^2 + (\rhoB)^2 < \frac{1}{\beta} \approx 52.632;
\end{equation}
i.e., NE will be unique for all physically relevant crowd densities.
However, this criterion becomes much more restrictive for higher (more realistic) values of $\beta$.  For example, at $\beta = 0.347$, the condition is now $(\rhoA)^2 + (\rhoB)^2 < 2.881$, which is definitely violated when $\rhoA + \rhoB \geq 2.41$ even if  $\max \left( \rhoA, \rhoB \right) \leq \hat{\rho} \approx 1.7$ and the velocity profiles are convex; see also Table \ref{table:RhoTable} and Figure \ref{fig:DensityRegions}.
Similarly, for the alternative disagreement penalty model based on $(\beta\bar{\rho})$,
the corresponding sufficient condition for the uniqueness of NE is
$\rhoA + \rhoB < \frac{1}{\beta} = \hat{\rho}.$  For $\beta = 0.347$, this means that the uniqueness might be lost when the total density of both crowds exceeds $2.881$.

Our numerical experiments show that the above conditions are actually sharp:
all densities from Region 2 in Figure \ref{fig:DensityRegions} seem to yield multiple NE.  However, the experiments also show that this only happens when $\bp$ and $(-\bq)$ are sufficiently close; e.g., at $\beta = 0.5$ in the $(\beta\bar{\rho}^2)$ model and $\rhoA = \rhoB = 1.2$, we only observe multiple NE when the angle between $\bp$ and $\bq$ is greater than 171 degrees.


\subsection{NE Examples.}
\label{ss:ExamplesNash}
\subsubsection*{Example 1}
Here we consider an example with $\bp = (1,0)$ and $\bq = (-1,0)$, with $\beta = 0.347$, $\rhoA = 1.68$, and $\rhoB = 0.72$
Numerically, we find three NEs at:
$$
\left(\begin{array}{c} a \\ b \end{array}\right)
 = \left(\begin{array}{c} 3.1416 \\ 0 \end{array}\right),
\left(\begin{array}{c} 3.0366 \\ 0.5208\end{array}\right),
\text{ and }   
\left(\begin{array}{c} 3.2466 \\  5.7623\end{array}\right).
$$

\begin{figure}[h]
	\centering
	\includegraphics[width=\textwidth]{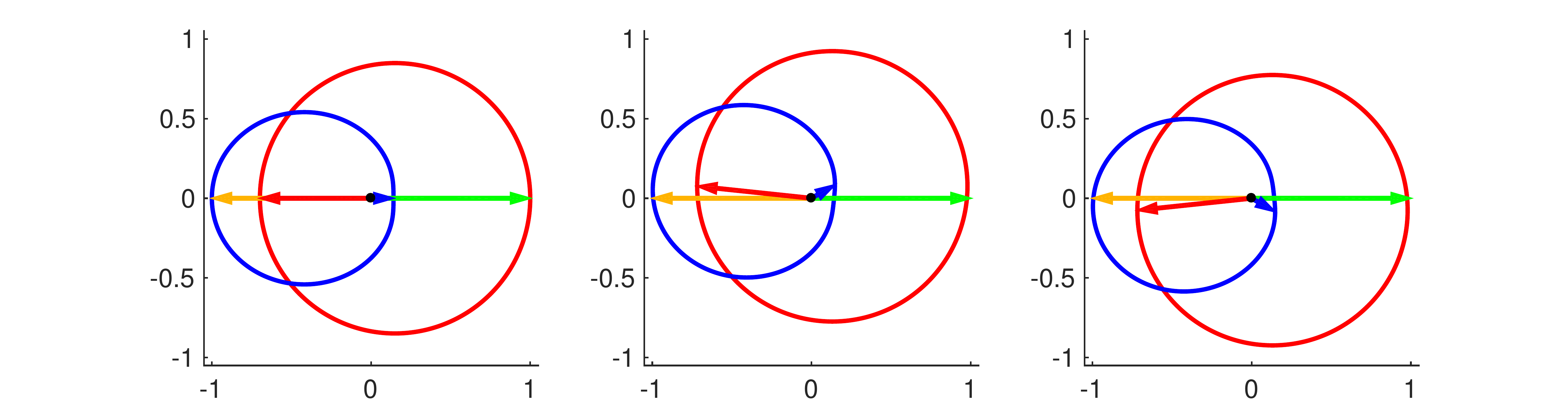}
	\caption{$-\bp$,$-\bq$, and velocity profiles for each NE in Example $1$.}
	\label{fig:NashEx1}
\end{figure}

A plot of $-\bp,-\bq$ and the velocity profiles for each NE can be found in Figure \ref{fig:NashEx1}. The orange arrow is $-\bp$ and the green arrow is $-\bq$. The velocity profile and choice of direction for crowd $A$ are plotted in red, while the same is plotted in blue for crowd $B$.

Example 1 is also summarized below in bimatrix form.  Both crowds try to maximize their respective outcomes based on \eqref{eq:game}; each entry contains the payoff for the ``row player'' (crowd A) and the ``column player'' (crowd B).  Of course, both of them have infinitely many options, but we include only those corresponding to Figure \ref{fig:NashEx1}.
The NE outcomes are on the diagonal, and their payoffs are different; so, the ``value'' of this game is not well-defined.
Moreover, since the crowds do not coordinate, they cannot know which NE control the other player will choose and this lack of coordination will frequently result in non-Nash (off-diagonal) outcomes.

\vspace{\baselineskip}

\begin{minipage}{0.9\textwidth}
\begin{center}
	\textbf{Nash Equilibria payoffs for Example 1}
	
	\begin{tabular}{c | c | c | c |}
		& 0 & 0.5208 & 5.7623 \\
		\hline
		3.1416 & 0.698, 0.141 & 0.715, 0.139 & 0.715, 0.139  \\
		\hline
		3.0366 & 0.695, 0.142 & 0.718, 0.147  & 0.705, 0.133 \\
		\hline
		3.2466 & 0.695, 0.142 & 0.705, 0.133 & 0.718, 0.147 \\
		\hline
	\end{tabular}
\end{center}
\end{minipage}

\vspace{\baselineskip}


\subsubsection*{Example 2}
\begin{figure}[h]
	\centering
	\includegraphics[width=\textwidth]{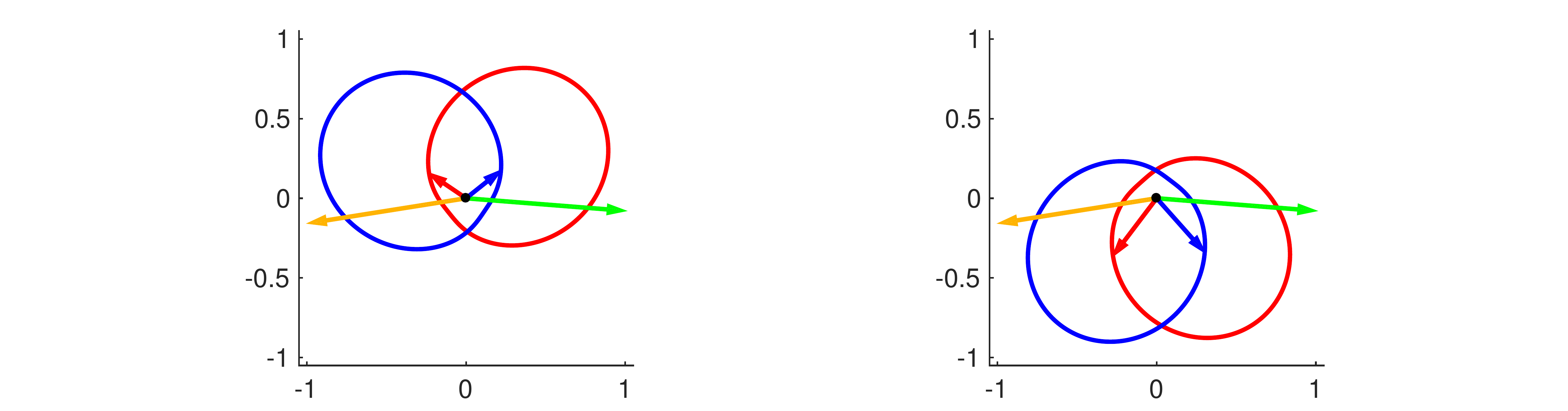}
	\caption{$-\bp$,$-\bq$, and velocity profiles for each NE in Example $2$.}
	\label{fig:NashEx2}
\end{figure}

Now for an example where $ \bp \not= -\bq$.
Here we take $\bp = (\cos(\pi/20),\sin(\pi/20))$ and $\bq = (\cos(39\pi/40),\sin(39\pi/40))$, with $\beta = 0.347$, $\rhoA = 1.68$, and $\rhoB = 1.68$.
Numerically, we now find two NEs, this time at:
$$
\left(\begin{array}{c} a \\ b \end{array}\right)
 = \left(\begin{array}{c} 2.5470 \\ 0.6732 \end{array}\right),
 \text{ and }
\left(\begin{array}{c} 4.0641 \\ 5.4393 \end{array}\right).
$$
A plot of $-\bp,-\bq$ and the velocity profiles for each NE can be found in Figure \ref{fig:NashEx2}

\vspace{\baselineskip}

\begin{minipage}{0.9\textwidth}
\begin{center}
	\textbf{NE payoffs for Example 2}
	
	\begin{tabular}{c | c | c |}
		& 0.6732 & 5.4393 \\
		\hline
		2.5470 & 0.205, 0.205 & 0.106, 0.105 \\
		\hline
		4.0641 & 0.105, 0.106 & 0.328, 0.328 \\
		\hline
	\end{tabular}
\end{center}
\end{minipage}

\vspace{\baselineskip}

From our numerical experiments, it seems that all cases with $\bp \not= -\bq$ and multiple NEs are qualitatively similar to Example 2. 
Since $\bp\not= -\bq$, we can think of $-\bp$ and $-\bq$ as splitting the plane into two sections, one where the angle between them is greater than $\pi$, and one where the angle is less than $\pi$. We always see two NEs, one in the section with angle greater than $\pi$, one in the other section. The NE on the side where the angle is less than $\pi$ then has higher payoffs for both crowds, dominating the other NE. 


\subsubsection*{Example 3}

In the previous examples, we showed that multiple NEs are possible even if the velocity profile is strictly convex and smooth. 
In this last example we aim to illustrate that a strictly convex but non-smooth velocity profile can result in infinitely many NEs.
Let us consider a speed function with ``disagreement penalty'' of the form:
$$f(\bar{\rho},\psi) = e^{-C\bar{\rho}\left(\pi^2-(\text{mod}(\psi,2\pi)-\pi)^2\right)}.$$

\begin{figure}[h]
\centerline{
$
\arraycolsep=0pt
\begin{array}{cccc}	
\includegraphics[width=.25\textwidth]{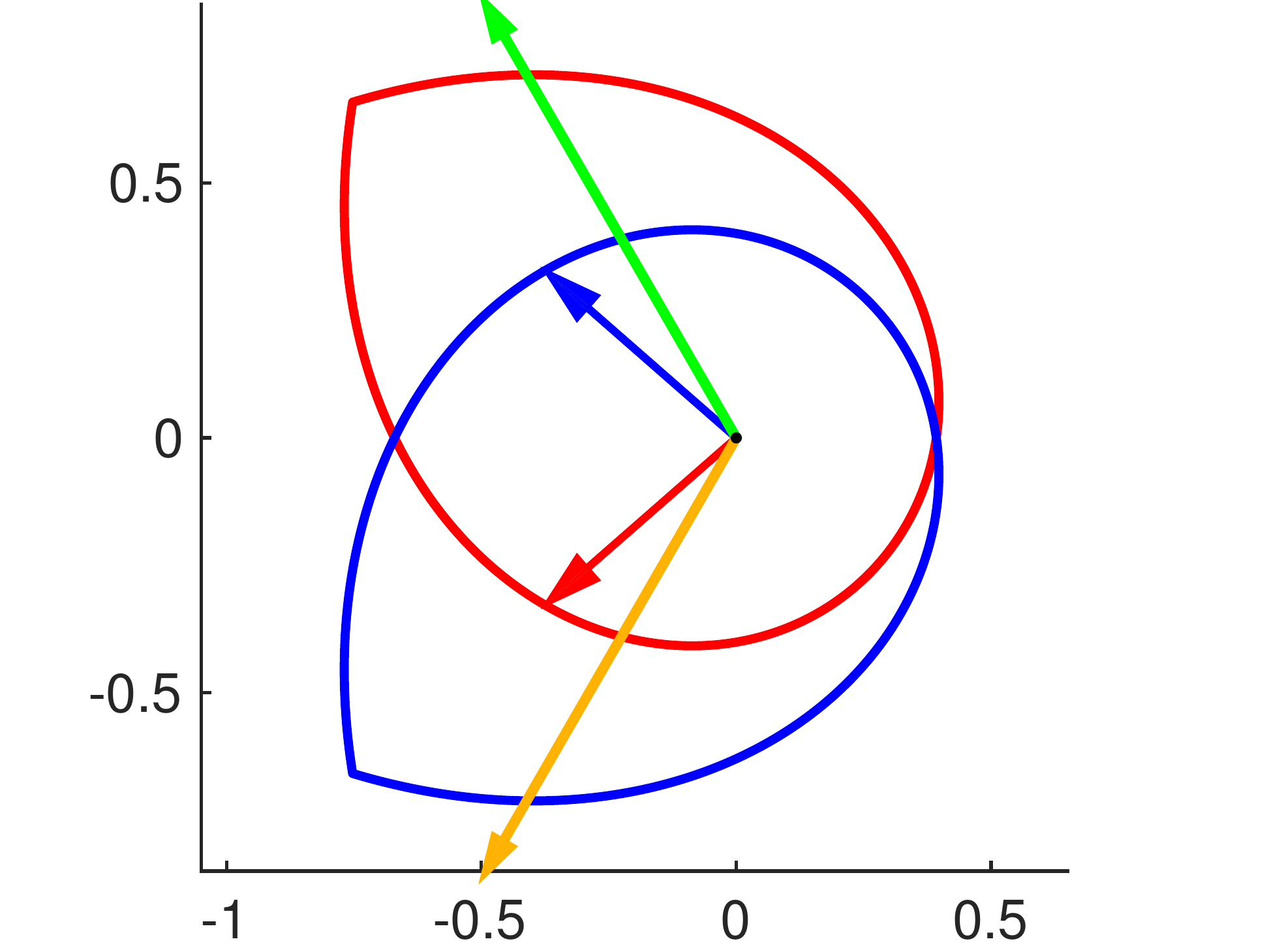} &
\includegraphics[width=.25\textwidth]{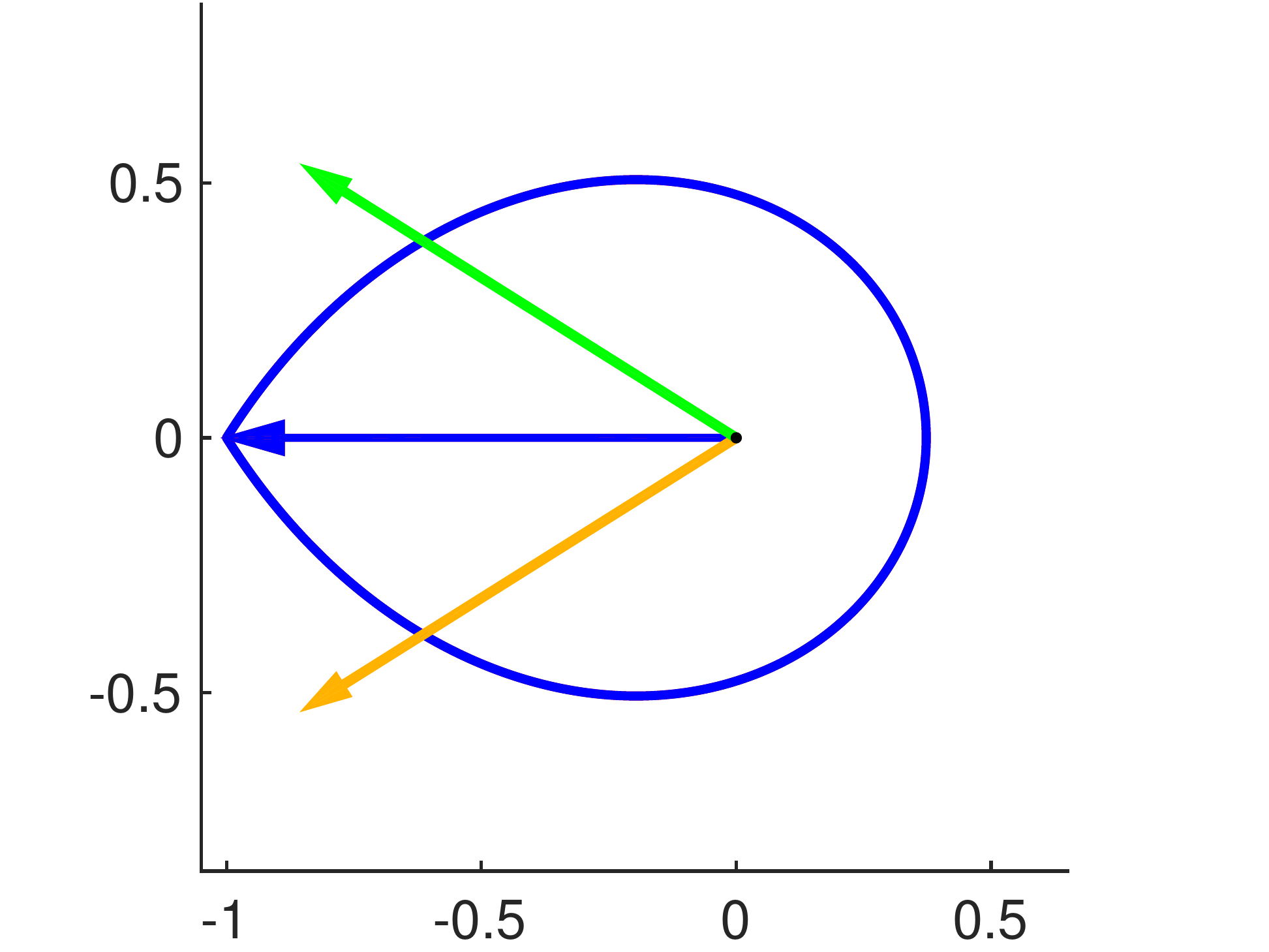} &
\includegraphics[width=.25\textwidth]{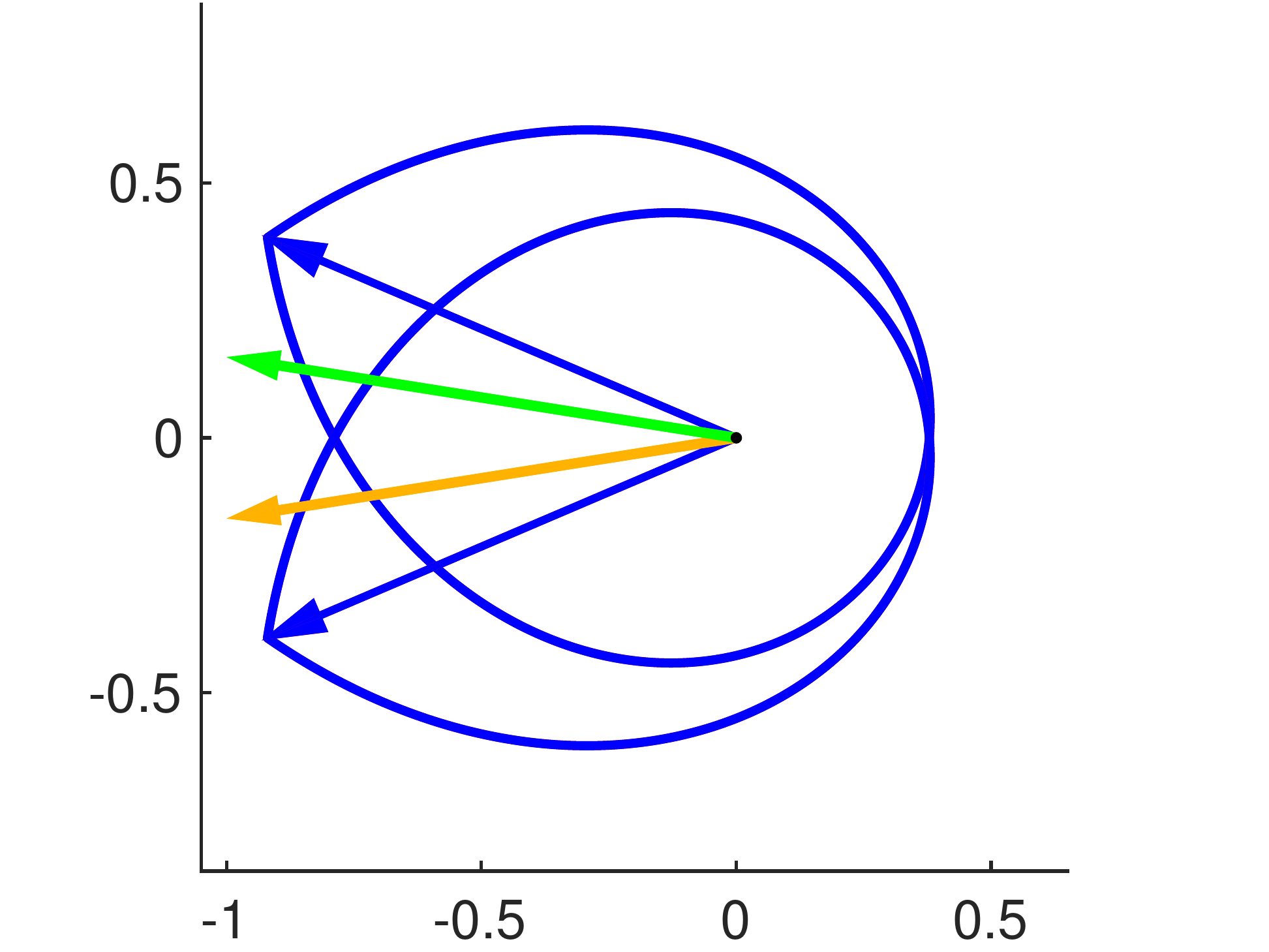} &
\includegraphics[width=.25\textwidth]{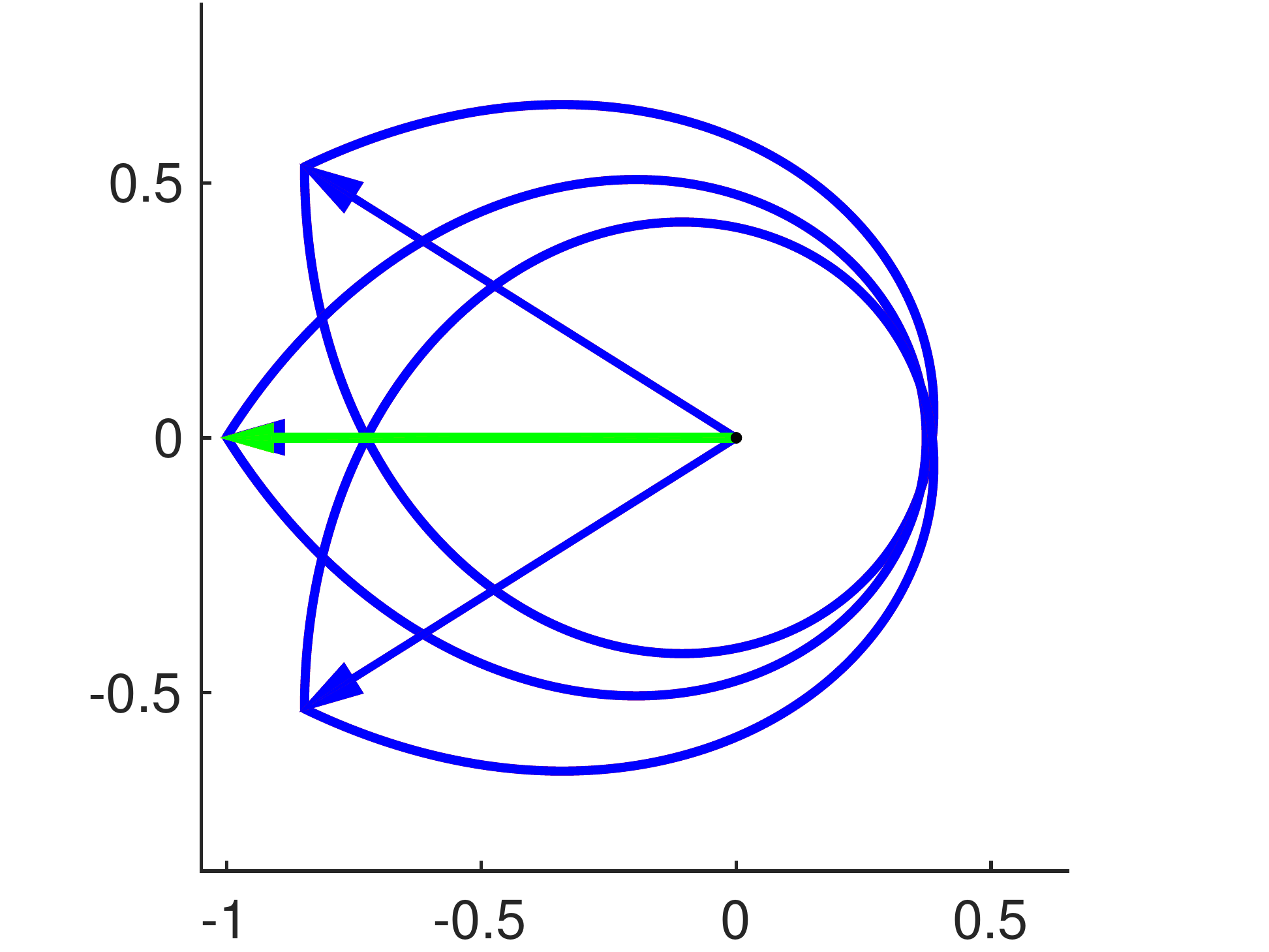}\\
(a) & (b) & (c) & (d)
\end{array}
$
}
\caption{Example 3: NEs for non-smooth velocity profiles as $(-\bp)$ in orange and $(-\bq)$ in green move closer together. Subfigures $(c)$ and $(d)$ show 
a few of infinitely many NEs.} 
\label{fig:NashEx3}
\end{figure}

The resulting velocity profiles form a strictly convex but non-smooth teardrop-like shape as seen in Figure \ref{fig:NashEx3}.
In the different subfigures of Figure \ref{fig:NashEx3}, we have consistently chosen $C = 0.1$ and $\rhoA = \rhoB = 1$, 
keeping $\bp$ and $\bq$ symmetric relative to the $x$-axis.
We explore how the NEs change as we decrease the angle $\theta$ between $(-\bp)$ and $(-\bq),$
which are shown by the  orange and green arrows respectively.

Figure \ref{fig:NashEx3}a shows the case $\theta = 2 \pi / 3$, with a unique NE, both crowds having distinct velocity profiles and different choices of optimal directions (shown in red and blue).
As  $\bp$ and $\bq$ move closer together, we still have a unique NE up until the critical angle $\theta_c,$ when the velocity profiles of both crowds (and the respective optimal directions) coincide; see Figure \ref{fig:NashEx3}b. 
A simple geometric argument shows that $\theta_c = \pi/2-\arctan(1/(2\pi C\bar{\rho})).$ 
For $\theta \in [0, \theta_c),$ the NE is no longer unique even though in each of them both velocity profiles coincide and both crowds choose the same ``optimal'' direction corresponding to the tip of the velocity profile.  Note that only some of these NEs are Pareto optimal (the ones where the tip lies in between $(-\bp)$ and $(-\bq)$), and many others are not; see the full range of Nash-profiles in Figure \ref{fig:NashEx3}c.
It is worth noting that the case $\bp=\bq$ can also be viewed as one crowd artificially divided into two, whose goals/destinations are actually the same.
In order for our model to remain consistent, we would hope that there is a unique NE corresponding to the optimal choice of direction for the isotropic (one-crowd) model.  Figure \ref{fig:NashEx3}d shows that this is actually {\bf not} the case and we have infinitely many NEs even with $\theta=0.$  (However, only one of them, corresponding to the vertically symmetric velocity profile and the ``isotropically optimal'' direction of motion $\frac{-\bp}{|\bp|}$, is now Pareto-optimal).
This suggests that the smoothness of the velocity profile is also a natural consistency condition even for the single crowd models.




\section{Numerical Implementation and Experiments.}
\label{s:NumericalImplementation}

The fully coupled system \eqref{eq:strobo_full} is typically discretized via explicit time-stepping: the current density $\rho$ is used to compute $\varphi$ in the current time slice, and the recovered $\bu_*$ is then used to find $\rho$ in the next time-slice via (\ref{eq:strobo_full}a, \ref{eq:strobo_full}d).
In \cite{Huang_2009} this approach was used with WENO-based discretizations for both the conservation law and the isotropic Eikonal PDE.
Our own simpler implementation relies on a Lax-Friedrichs scheme \cite{lax1954weak} for the hyperbolic conservation law (\ref{eq:strobo_full}a) and a first-order accurate semi-Lagrangian discretization of the HJB equation (\ref{eq:strobo_full}b).

Modeling the choice of optimal directions presents two main challenges.
The first is a choice of suitable discretization for the anisotropic / static HJ PDEs, yielding a large system of coupled (discretized) equations in each time slice.
The second is an efficient method for solving this discretized system.  In the following discussion we primarily focus on the latter challenge and also emphasize the connection to our consistency criteria covered in previous sections.  We will first discuss the numerical methods for a single HJB PDE (\ref{eq:strobo_full}b) in Section \ref{ss:numerics},
followed by a discussion of numerics for two coupled HJI equations in Section \ref{ss:NashNumerics}, and two-crowd simulation results in Section \ref{ss:NumericalExperiments}.

\subsection{Numerical solution of HJB for a single anisotropic crowd.}
\label{ss:numerics}
If $\rho$ is known and $t>0$ is fixed, then $\varphi$ is the viscosity solution of the Hamilton-Jacobi-Bellman equation (\ref{eq:strobo_full}b), and $\bu_*$ can be recovered from (\ref{eq:strobo_full}c).
In the special case of isotropic dynamics (i.e. when $\bv(\rho(\x,t),\bu) = v(\rho(\x,t))\bu$ and $\bu \in S^1$), this HJB PDE reduces to an Eikonal equation \eqref{eq:Eikonal}, for which many competing efficient algorithms have been developed in the last 20 years.
Dijkstra-like (``Fast Marching'') methods \cite{tstitsiklis1996, sethian1996fast, SethBook2} provide non-iterative solvers by dynamically recovering the correct ordering of discretized equations and solving them one at a time.  These methods have been extended to meshes on manifolds \cite{KimmSethTria, SethVlad1}, higher-order accurate Eulerian \cite{SethSIAM, Renzi3} and semi-Lagrangian discretizations \cite{CristianiFalcone_2007}, and more general anisotropic cases \cite{sethian2003ordered, AltonMitchell2, Mirebeau3}.
An alternative (``Fast Sweeping'') approach is based on alternating through a set of natural geometric orderings to speed up the convergence of Gauss-Seidel iterations \cite{BoueDupuis, zhao2005fast, TsaiChengOsherZhao, KaoOsherQian, Li_DG, Zhang_Eik, Zhang_HJB, BenamouLuoZhao}.  Each of these two approaches is particularly efficient on its own subset of problems.   More recently, hybrid (two-scale) methods  have been introduced in \cite{ChacVlad1,ChacVlad2} to combine the advantages of both marching and sweeping.  These papers also include a review of several other ``fast Eikonal solvers'' (mirroring  the logic of label-correcting algorithms on graphs) and their parallel implementations.

Our own simulations use Fast Sweeping -- primarily because this approach is easier to extend to multi-crowd scenarios discussed in section \ref{ss:NashNumerics}.

In principle, the computational cost of solving (\ref{eq:strobo_full}b) in each time-slice can be further reduced by only computing $\varphi$ on the relevant subset of $\cdomain$.  Such ``causal domain restriction'' techniques are reminiscent of the classical A* algorithm on graphs \cite{Hart_Astar} and have been also extended to Eikonal equations; see \cite{clawson2014causal} and references therein.  In pedestrian flow modeling, similar ideas have been used to speed up the ``dynamic floor field'' computation in \cite{Hartmann_2014}.

For the general/anisotropic version of the pedestrian direction field problem, it usually is not possible to write down the analytic expression for $\bu_*$ even if $\nabla \varphi$ is already known.  The HJB equation (\ref{eq:strobo_full}b) is then most naturally discretized in a semi-Lagrangian framework \cite{Falcone_book}, with $\bu_*$ recovered at each gridpoint through a numerical optimization over all $\bu \in \mathcal{U}$.

This numerical optimization often presents a significant computational bottleneck, particularly so for non-local models, in which evaluating the velocity $\bv$ for any specific $\bu$ is already expensive.
We recall that in \cite{Cristiani_2015} the velocity is computed as the sum of a chosen velocity $\bu$ and an interactional velocity $\bw_i$; see formula \eqref{eq:nonLocal_velocity}.
Computing the value of $\bw_i$ requires calculating an integral over a sensory region $S(\x,\bu)$ that changes with each choice of $\bu$.
However, if the density is approximately linear within the sensory region $S(\x,\bu)$, our work in Appendix \ref{ss:Appendix1} shows how to evaluate this integral analytically, which could be used to speed up the numerics.

Aside from being computationally expensive, the numerical optimization to find $\bu_*$ may ``mask'' some of the possible model inconsistencies: if the argmax in (\ref{eq:strobo_full}c) is not unique (e.g., when the velocity profile is not strictly convex), this can easily go unnoticed and the simulation results will become dependent on the specific implementation of optimization on $\mathcal{U}.$  To avoid these issues, it is important to verify the convexity of the speed profile either a priori or at run-time (based on the crowd densities observed in the simulation).

\subsection{Implementation notes for two-crowd models.}
\label{ss:NashNumerics}

In Section \ref{s:Nash}, we discussed uniqueness of Nash Equilibria for the system of coupled Hamilton-Jacobi-Isaacs PDEs \eqref{eq:HJI_1} and \eqref{eq:HJI_2}.
For numerical simulations, the semi-Lagrangian discretization of these PDEs on the grid results in a pair of coupled optimization problems at each gridpoint $(i,j)$:
\begin{equation}
\varphiA_{ij} \; = \; \min_{\buA_{ij} \in S^1} G\left(\rhoA_{ij}, \rhoB_{ij}, \buB_{ij}, N \varphiA_{ij} \right),
\label{eq:discHJI_A}
\end{equation}
\begin{equation}
\varphiB_{ij} \; = \; \min_{\buB_{ij} \in S^1} G\left(\rhoA_{ij}, \rhoB_{ij}, \buA_{ij}, N \varphiB_{ij} \right),
\label{eq:discHJI_B}
\end{equation}
where $N \varphiA_{ij}$ refers to the values of $\varphiA$ at the neighboring gridpoints to $(i,j)$.
Ideally, we would like to know that the system (\ref{eq:discHJI_A}-\ref{eq:discHJI_B}) has a unique Nash Equilibrium regardless of values of
$\rhoA_{ij}, \rhoB_{ij}, N \varphiA_{ij},$ and $N \varphiB_{ij}$ whenever the same is true for (\ref{eq:HJI_1}-\ref{eq:HJI_2}) regardless of the vectors $\nabla \varphiA$ and $\nabla \varphiB$.
Our numerical experiments confirm this equivalence, but it has not been proven rigorously so far.

These coupled discretized equations have to be solved at every gridpoint.  A natural approach for doing this is to mimic the structure of the ``best reply'' map in Section \ref{s:Nash}.  If $\varphi$ for one of the crowds is assumed to be known, then $\varphi$ for the other crowd can be solved as in the one-crowd case (e.g., by Fast Sweeping; see Section \ref{ss:numerics}).
So, we can alternate through solving equations \eqref{eq:discHJI_A} and \eqref{eq:discHJI_B}, freezing the current version of $\varphiA$ and $\buA$ to recover $\varphiB$ and $\buB$,
which are then used to update $\varphiA$ and $\buA$, and so on.
When the (discretized) best reply map is a contraction, its fixed point is the unique Nash equilibrium and this iterative process converges to it.
The examples in Section \ref{ss:NumericalExperiments} show simulations of two crowds computed in this manner which exhibit qualitative behavior similar to real crowds.
If we initialize the discretized best reply map with some initial choice of directions (e.g. from the previous time step), then the convergence of this map is proof that both crowds are in fact using Nash equilibrium strategies.
However, this approach is not sufficient to guarantee that this Nash equilibrium is also unique, even if the simulation appears plausible from a phenomenological perspective.

A more careful approach is to iterate over the gridpoints in the outer loop,  freezing $\left( N \varphiA_{ij}, N \varphiB_{ij} \right)$ and updating $(\varphiA_{ij}, \varphiB_{ij})$ simultaneously.
We can discretize the control space $S^1$ (considering a fixed/finite number of directions of motion) and then find the Nash equilibria for the resulting finite bi-matrix game.
One advantage is that we can directly verify the number of Nash equilibria, though it might also change as we refine the discretization of $S^1.$

Since the coupled system (\ref{eq:HJI_1}-\ref{eq:HJI_2}) has to be solved in every time-slice, it is desirable to speed up the above procedures as much as possible.
One heuristic approach is based on updating the crowd direction fields by measuring the angle $\Psi$ between our {\bf currently} chosen direction and the direction used by the other crowd in the {\bf previous} time-slice.
This idea (similar to the approach used in \cite{Jiang_2009} and the subsequent papers) effectively decouples the PDEs, speeding up the numerics and side-stepping all Nash-equilibria issues.
However, there is no proof that the resulting evolution of crowd densities is the same as one would recover from actually solving the coupled system (\ref{eq:HJI_1}-\ref{eq:HJI_2}).

\subsection{Numerical Experiments with Two Crowds.}
\label{ss:NumericalExperiments}
This subsection contains a few simulation results based on our 
implementation of the anisotropic two-crowd model found in \cite{Jiang_2009}, discussed previously in Sections \ref{ss:SpeedProfiles_Inter} and \ref{ss:Nash_Inter}.
As a reminder, the corresponding speed functions are:
\begin{equation*}
\vA \left(\rhoA,\rhoB,\buA,\buB \right)  = \bar{v} e^{-\alpha \left(\rhoA+\rhoB \right)^2} f(\rhoB, \psi)
\; \text{ and }  \;
\vB \left(\rhoA,\rhoB,\buA,\buB \right)  = \bar{v} e^{-\alpha \left(\rhoA+\rhoB \right)^2} f(\rhoA, \psi)
\end{equation*}
where $f$ is given by:
\begin{equation*}
f(\bar{\rho}, \psi)  = e^{-\beta \left(1-\cos \psi \right)\bar{\rho}^2}.
\end{equation*}
In all examples, we will use $\alpha = 0.075$ and $\bar{v} = 1$ meter/second on a $100$m by $100$m square domain $\Omega$.
Our code also verifies that the models stay consistent throughout the simulations; i.e., the densities stay in Region 1 defined in Figure \ref{fig:DensityRegions}a.

\subsubsection*{Example 1: Difference in Direction Fields}
\begin{figure}[h]
\centerline{
$
\arraycolsep=-4pt
\begin{array}{ccccc}	
\includegraphics[width=.095\textwidth]{Figures/ColorBar.pdf} &
\includegraphics[width=.28\textwidth]{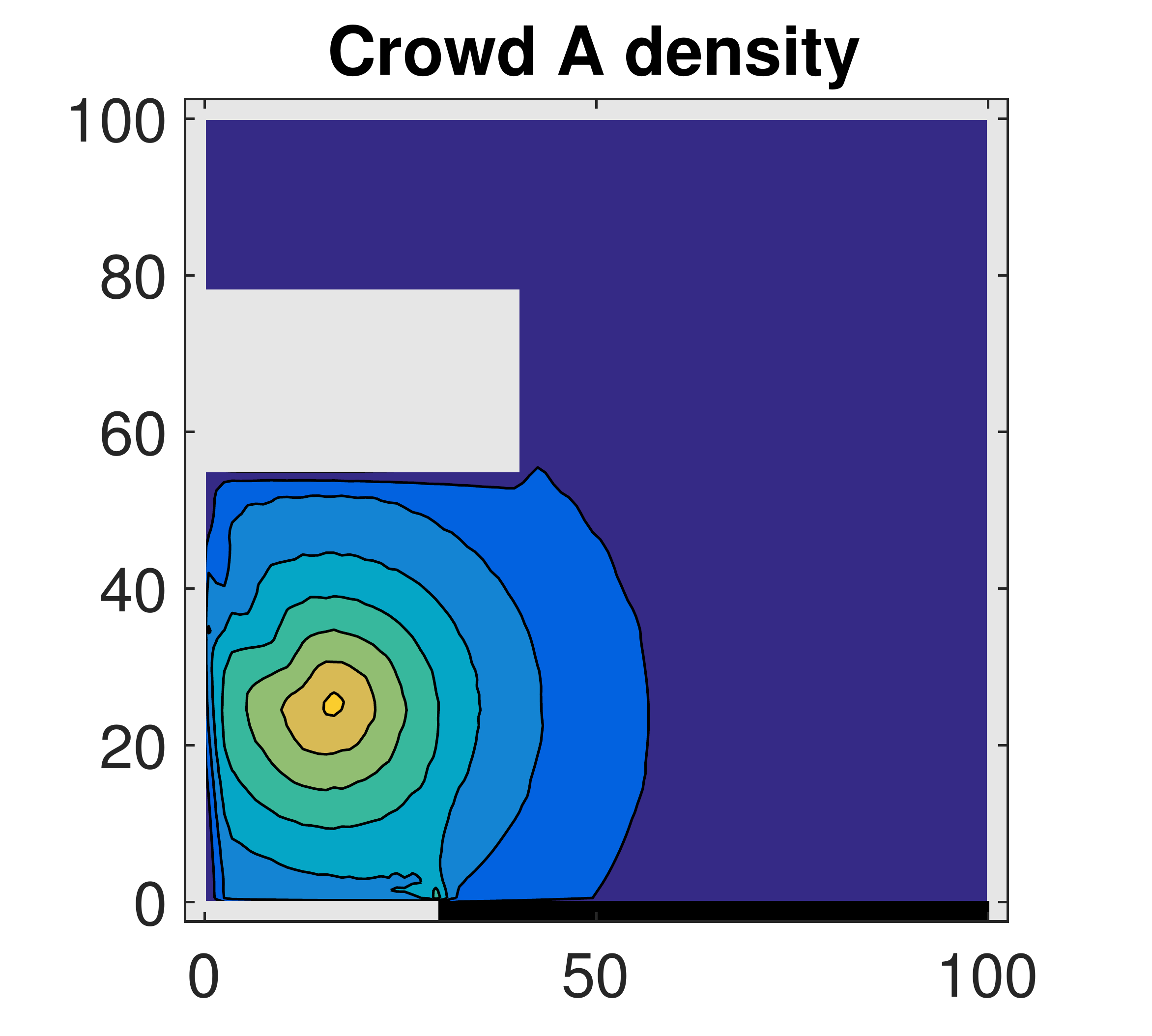} &
\includegraphics[width=.28\textwidth]{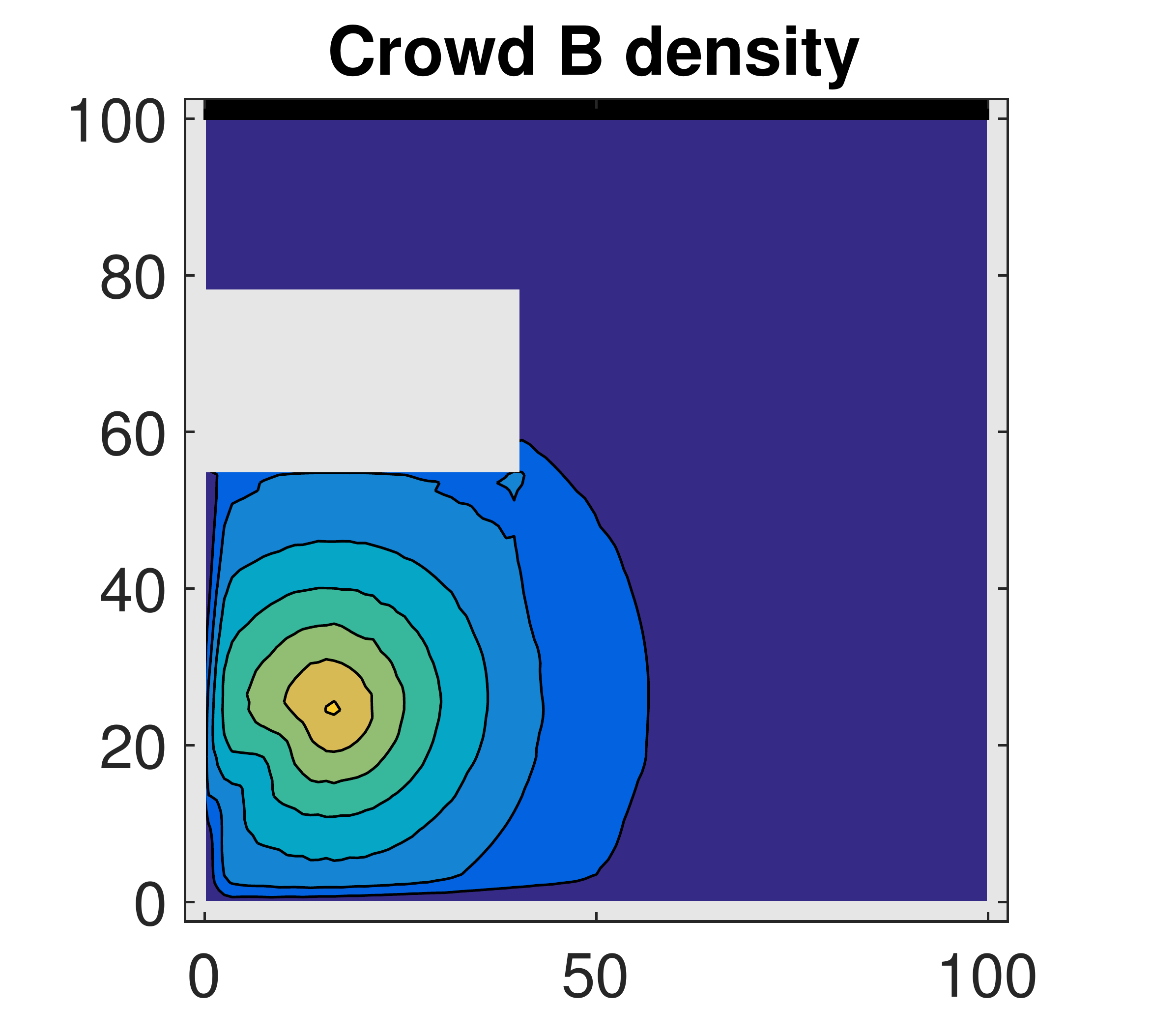} &
\includegraphics[width=.28\textwidth]{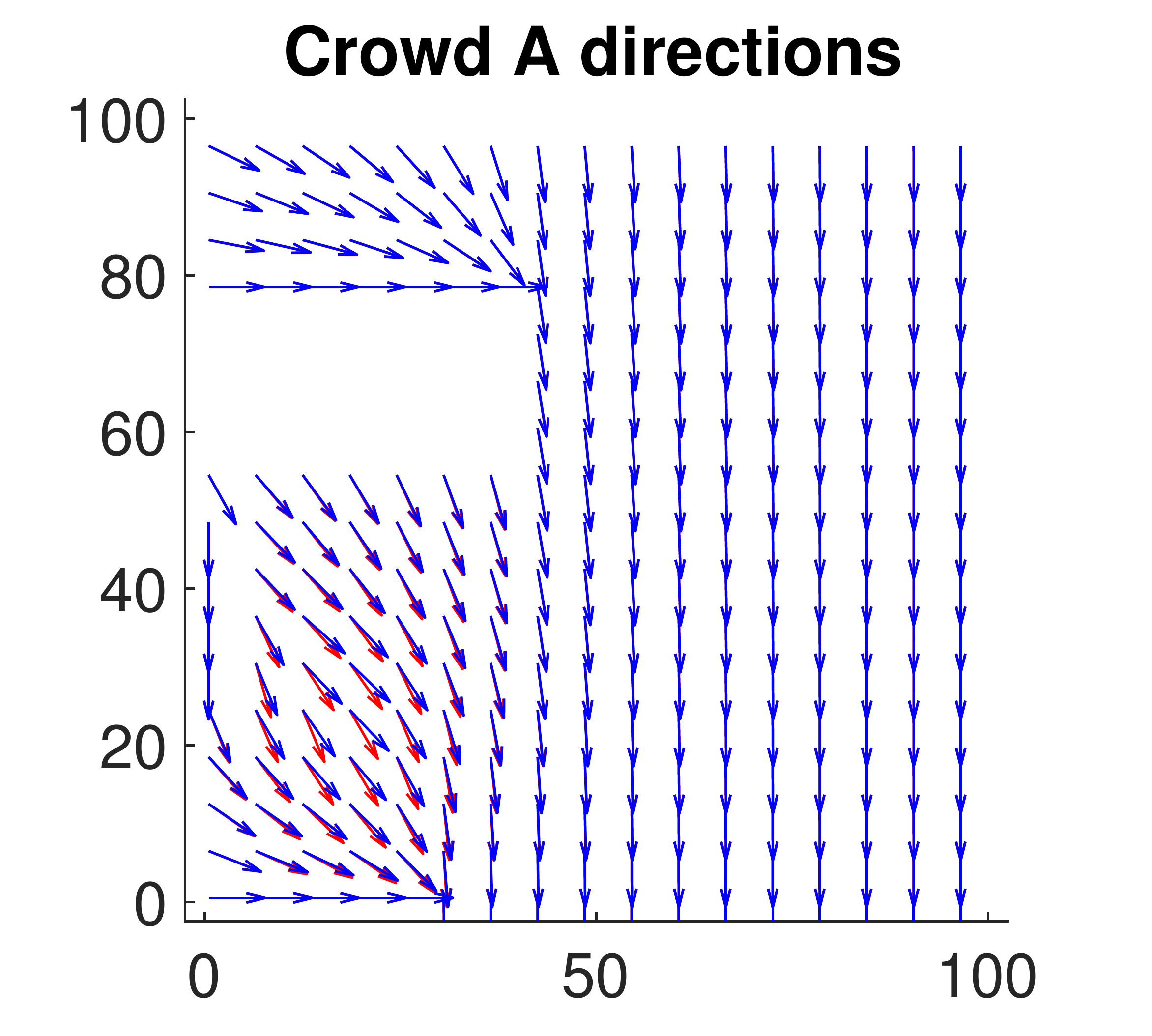} &
\includegraphics[width=.28\textwidth]{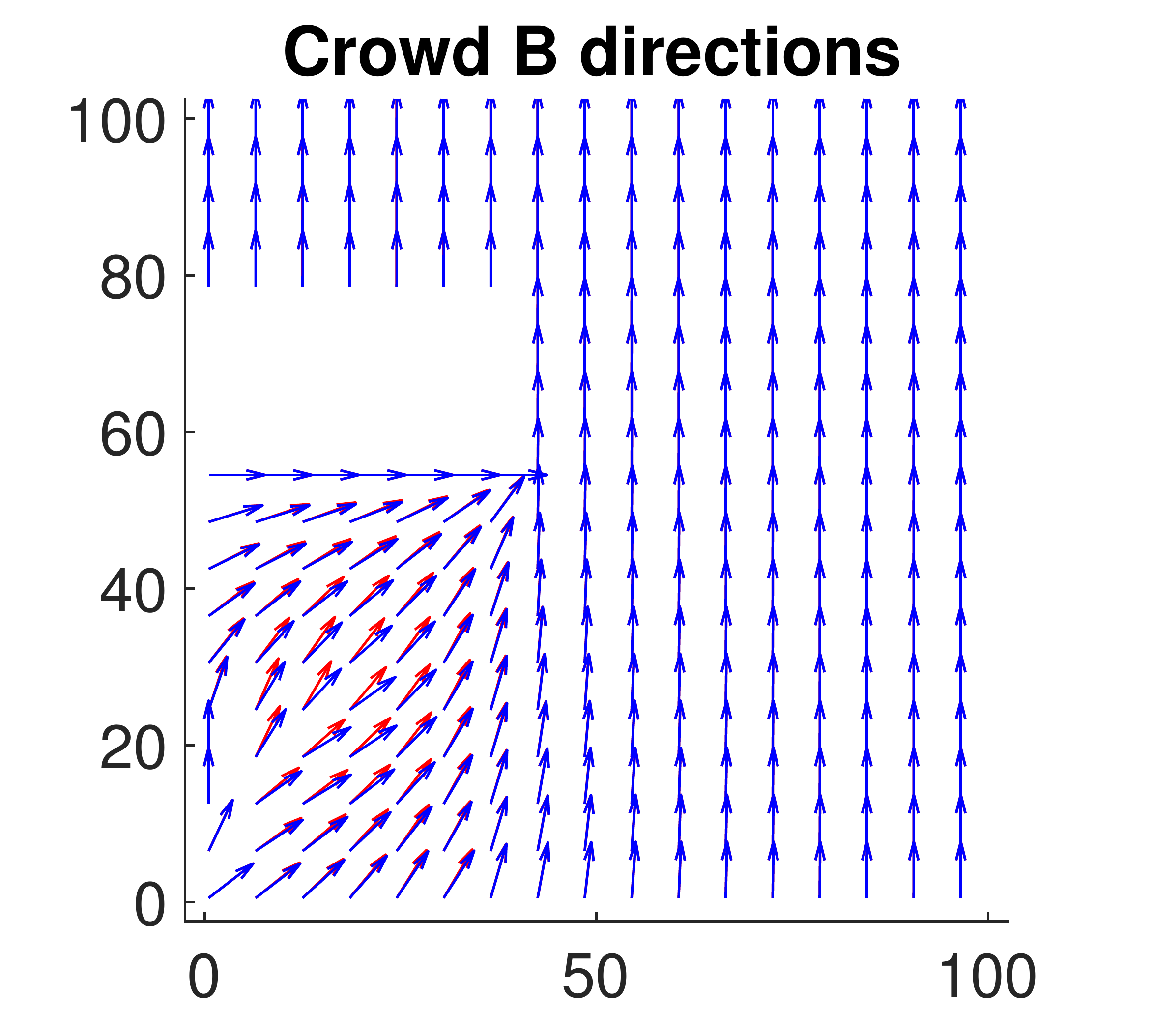} \\
& (a) & (b) & (c) & (d)
\end{array}
$
}
\caption{Two interacting crowds aiming for different exits on a domain with a single obstacle.
Subfigures (a) and (b): contour plots of $\rhoA$ and $\rhoB$ at $t=1.5$ with the respective exits shown in black.
Subfigures (c) and (d): their corresponding optimal direction fields (in blue) and gradient fields (in red).
Computed for $\beta = 0.347$ on a $100$x$100$ grid with timestep $\Delta t = 0.5.$
}
\label{fig:diff}
\end{figure}

This example is meant to highlight the difference between optimal and gradient path-planning.
In Figure \ref{fig:diff}, crowds A and B are trying to reach exits on the South and North sides of the domain, respectively.
The crowds are initially intermixed and their optimal direction fields (shown in blue) are different from the gradient fields of
their respective value functions (shown in red).
The red directions appear to indicate shorter paths toward the desired exits, but using them would be suboptimal since the corresponding ``disagreement penalty'' to the speeds would be higher.
This can be interpreted as a more ``organic'' version of the example 
in Figure \ref{fig:SwimmerRiver}.

\subsubsection*{Exmaple 2: Lane Formation}
\begin{figure}[h]
\centerline{
$
\arraycolsep=-4pt
\begin{array}{ccccc}	
\includegraphics[width=.095\textwidth]{Figures/ColorBar.pdf} &
\includegraphics[width=.28\textwidth]{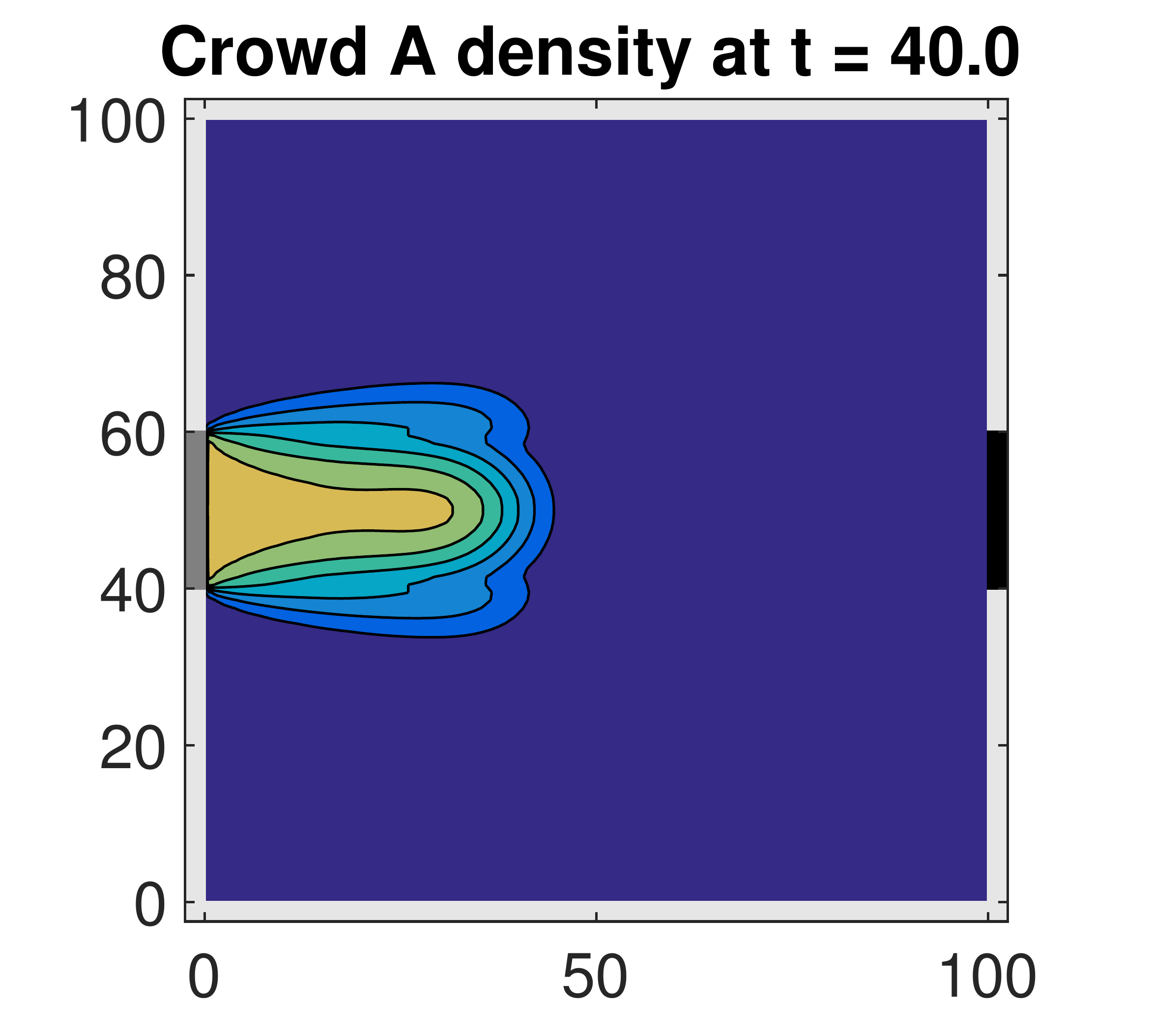} &
\includegraphics[width=.28\textwidth]{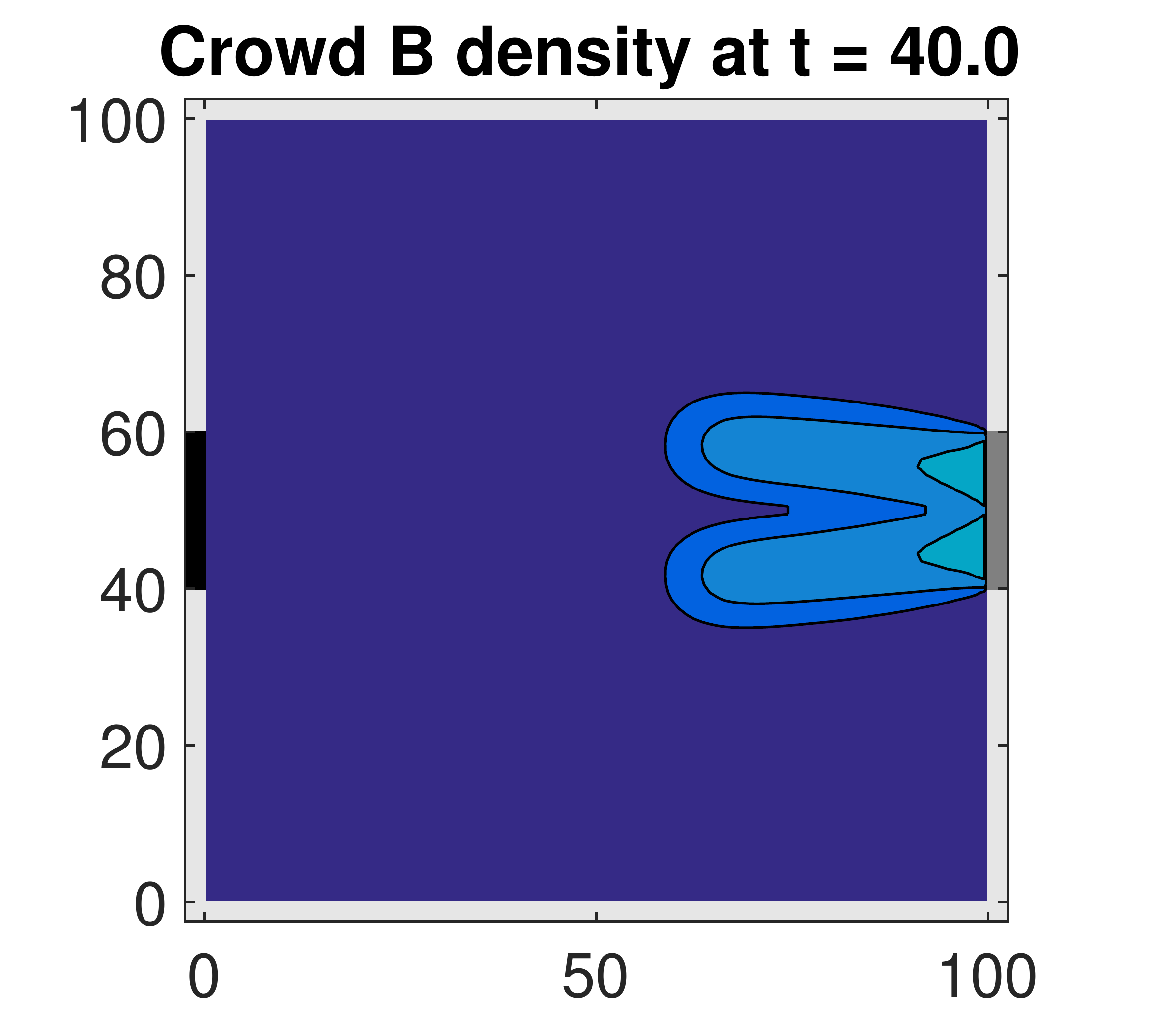} &
\hspace*{5mm}
\includegraphics[width=.28\textwidth]{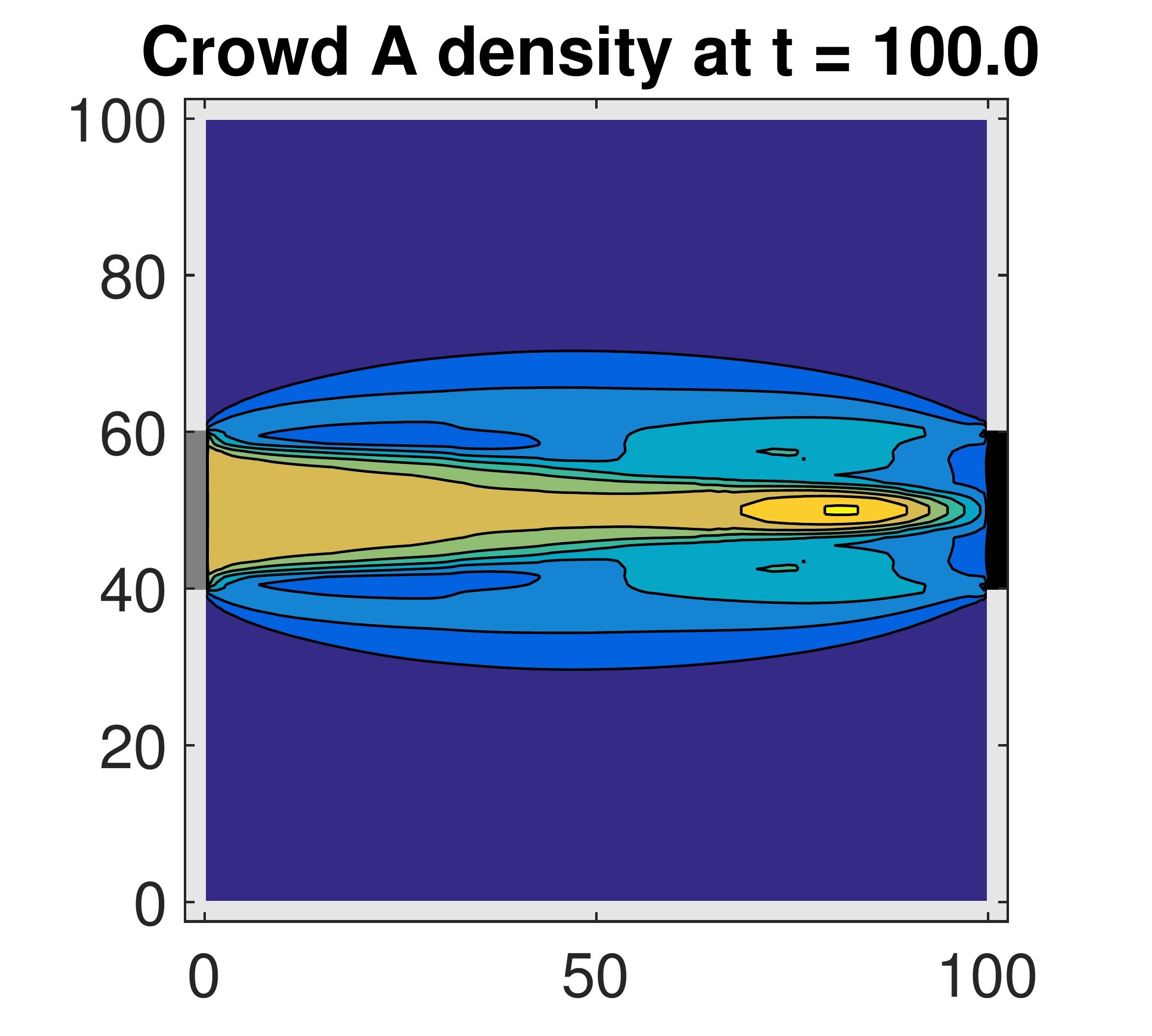} &
\includegraphics[width=.28\textwidth]{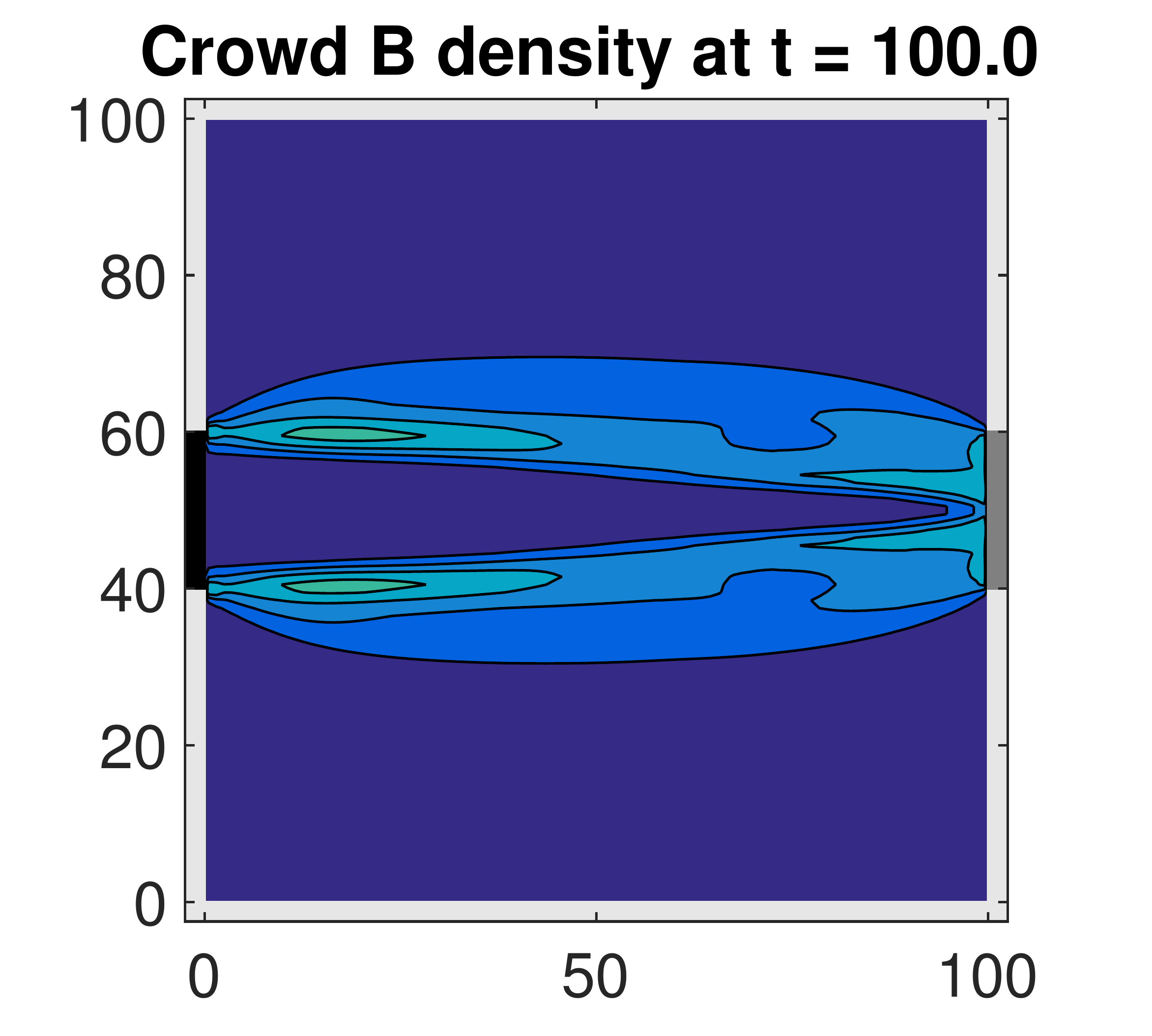} \\
& (a) & (b) & (c) & (d)
\end{array}
$
}
\caption{A simulation of two crowds, taken at time slices of $t=40$ and $t=100$. Crowd A enters from the left side of the domain (shown in dark gray) and is trying to reach an exit on the right (shown in black), while Crowd B enters from the right side of the domain and is trying to reach an exit on the left. Computed for $\beta=0.178$ on a $100$x$100$ grid with a timestep of $\Delta t = 0.5.$}
\label{fig:LaneFormation}
\end{figure}

This example is meant to illustrate lane formation of two crowds in an intercrowd model based on \cite{Jiang_2009}.
Each crowd flows in from an entrance on the side of the domain (shown in dark gray), with Crowd $A$ flowing in from the left at a rate of 0.8 $\text{ped}\cdot m^{-1} \cdot s^{-1}$ and Crowd $B$ flowing in from the right at a rate of 0.4  $\text{ped}\cdot m^{-1} \cdot s^{-1}$.
Each crowd then aims to reach a target on the opposite side of the domain (shown in black).
As can be seen in Figure \ref{fig:LaneFormation}, both crowds form lanes to avoid each other.  Since this is similar to behavior
of real crowds, such phenomena are often used to argue for the validity of PDE models.
It is worth noting that the lanes do not separate completely, so correctly resolving the anisotropic interactions here is quite important.
We can quantify this effect by normalizing the two densities and calculating the {\em overlapping coefficient} OVL \cite{inman1989overlapping} as follows:
\begin{equation*}
\text{OVL}(t) = \int_{\Omega} \min \left( \frac{\rhoA (\x,t)}{\int_{\Omega} \rhoA(\x,t)d\x}, \frac{\rhoB (\x,t)}{\int_{\Omega} \rhoB (\x,t)d\x} \right) d\x
\end{equation*}
This overlapping coefficient measures the similarity of the two densities, with $\text{OVL}=1$ if $\rhoA = \rhoB$ and $\text{OVL}=0$ if the supports of $\rhoA$ and $\rhoB$ are disjoint.
For the densities in subfigures \ref{fig:LaneFormation}(c-d), the overlapping coefficient is $\text{OVL}(100) = 0.5841$.
While some of this overlap is due to numerical viscosity from solving Equation \eqref{eq:strobo_full}(a) via Lax-Friedrichs, experimentally we see that OVL does not change significantly under grid refinement.

\subsubsection*{Example 3: Intersection}
\begin{figure}[h]
\centerline{
$
\arraycolsep=-6pt
\begin{array}{ccccc}	
\includegraphics[width=.095\textwidth]{Figures/ColorBar.pdf} &
\includegraphics[width=.28\textwidth]{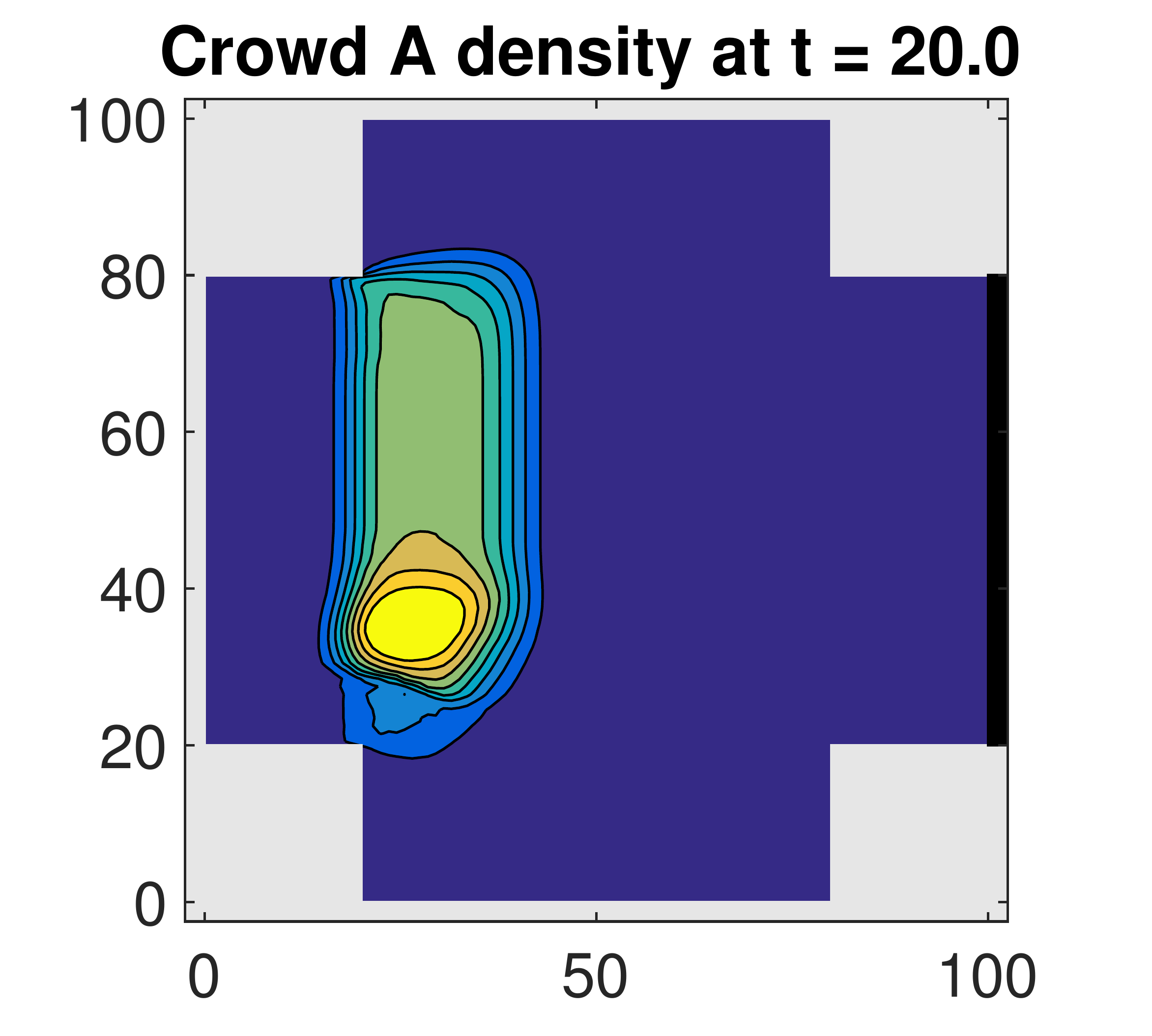} &
\includegraphics[width=.28\textwidth]{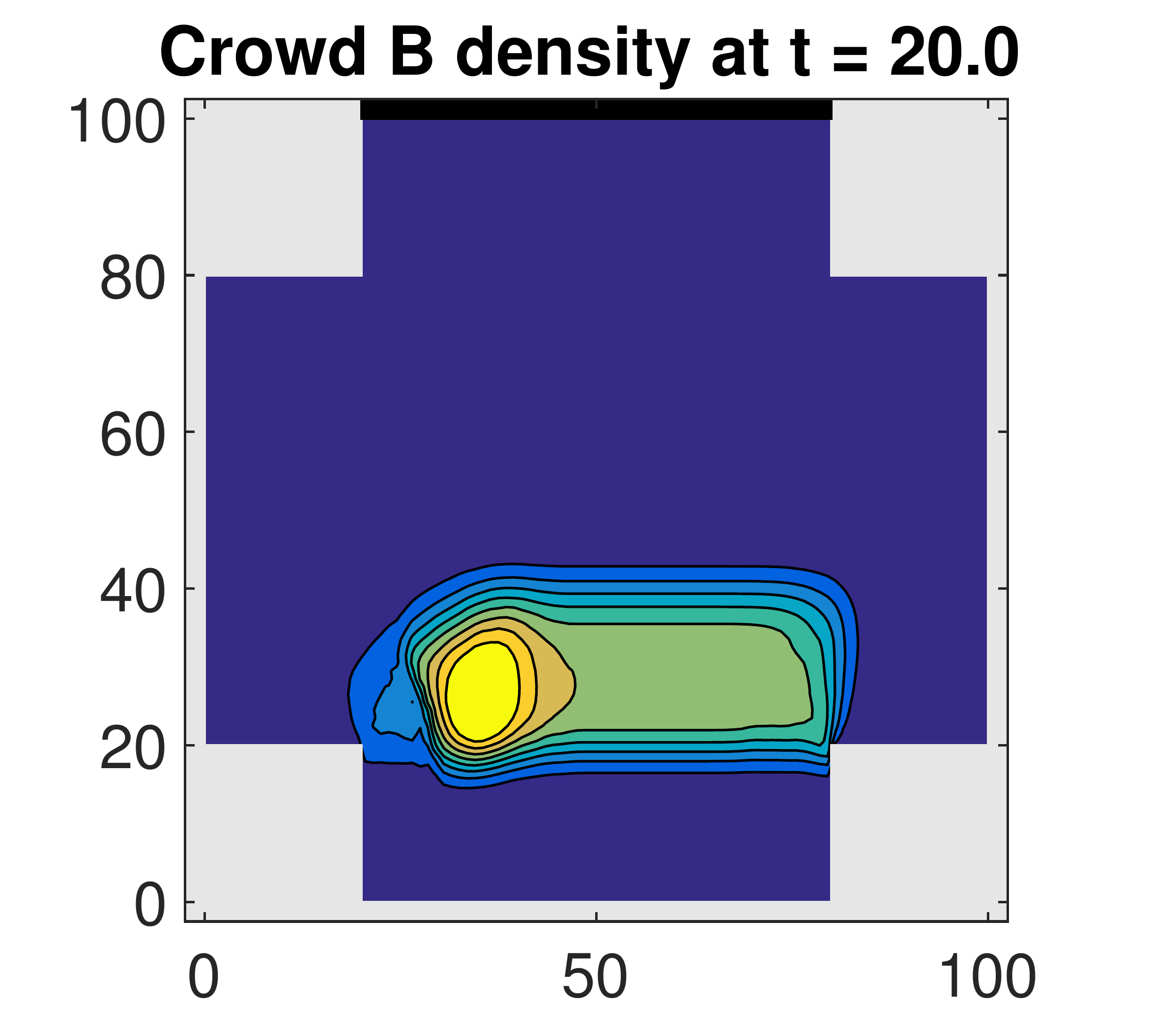} &
\hspace*{5mm}
\includegraphics[width=.28\textwidth]{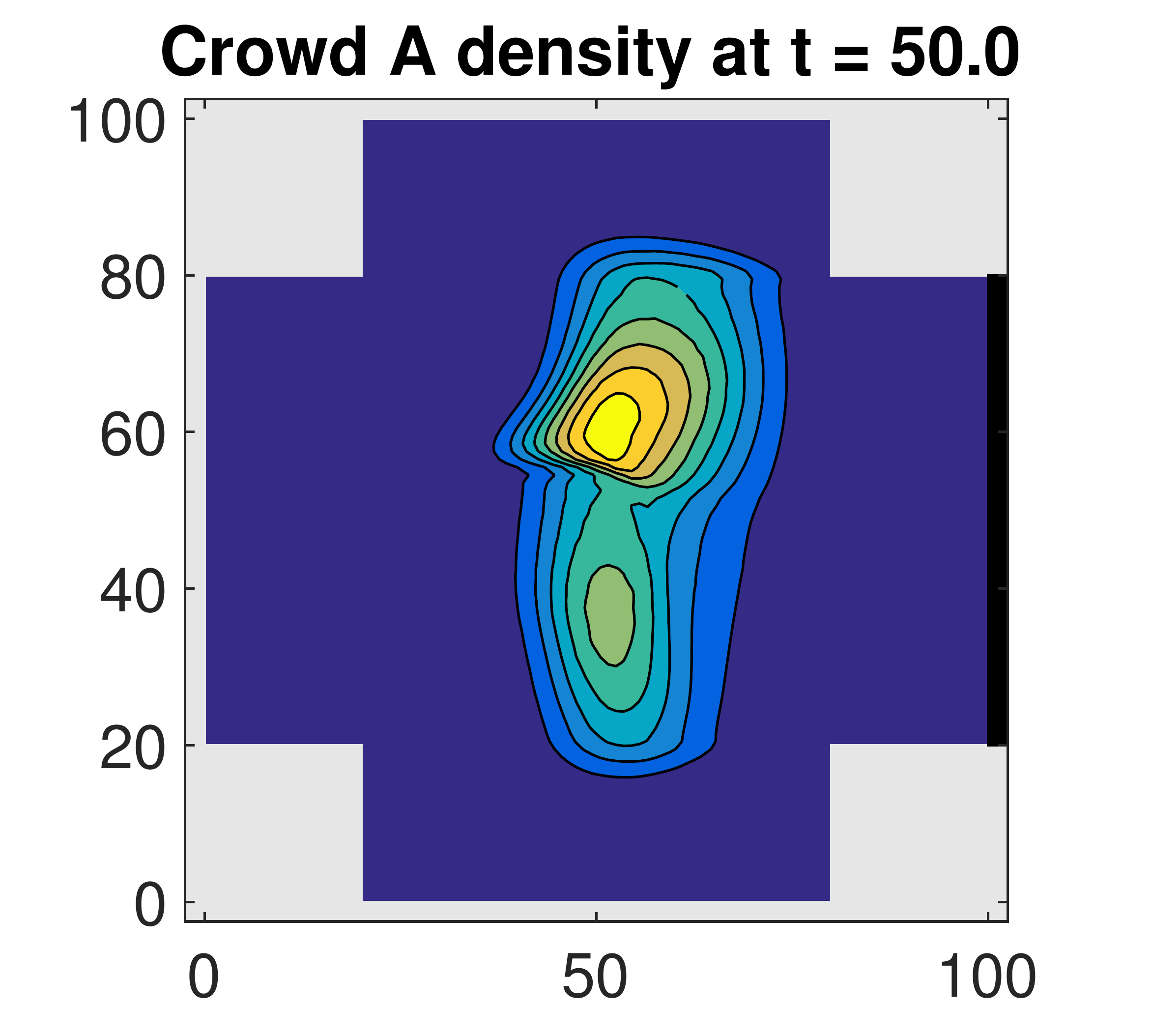} &
\includegraphics[width=.28\textwidth]{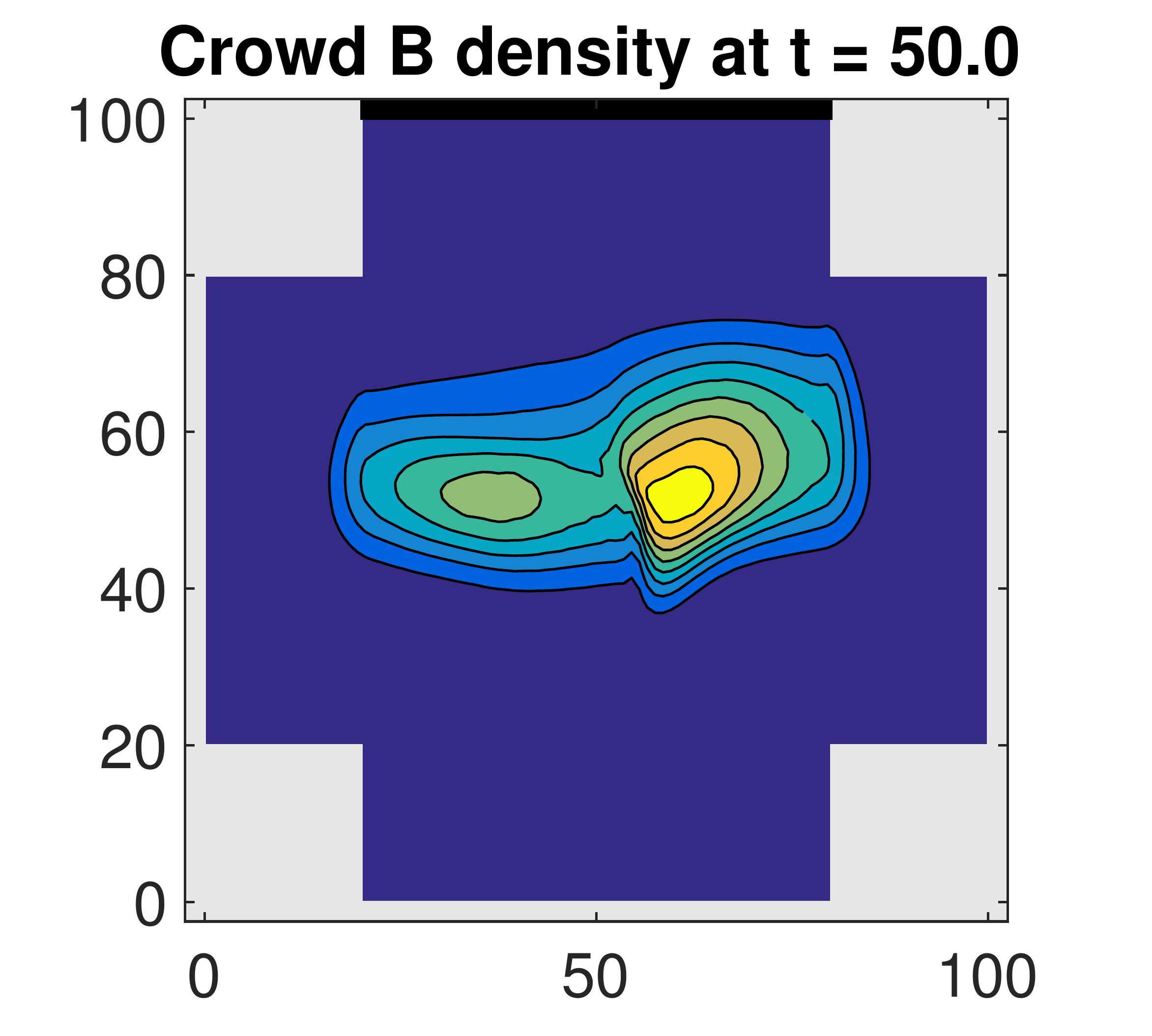} \\
& (a) & (b) & (c) & (d)
\end{array}
$
}
\caption{Two crowds meeting at an intersection (at $t=10$ and $t=50$). Crowd A begins on the left side of the domain and is trying to reach an exit on the right side of the domain (shown in black), while Crowd B begins at the bottom of the domain and is trying to reach an exit at the top. Computed for $\beta=0.178$ on a $100$x$100$ grid with a timestep of $\Delta t = 0.5.$}
\label{fig:Intersection}
\end{figure}

As in Example 2, each crowd begins on one side of the domain and aims to reach a target at the opposite side. 
However, instead of moving in opposite directions to each other (as in Example 2), they now meet at a right-angle intersection in the middle of the domain (seen in Figure \ref{fig:Intersection}).

In subfigures \ref{fig:Intersection}(c-d) the crowds now again overlap significantly, with an overlapping coefficient of $\text{OVL}(50) = 0.4538$.
Moreover, the most of this overlap happens near the center of $\Omega,$ where the difference in direction fields between the gradient and optimal path planning is the largest,
making it important to resolve the anisotropic interactions correctly.

\section{Conclusions.}
\label{s:Conclude}
We have reviewed several anisotropic crowd flow models, highlighting the geometric interpretation of ``optimality'' in choosing the pedestrians' direction field.  We have derived internal consistency criteria to guarantee that this direction field is uniquely defined almost everywhere on $\domain$ and showed that the strict convexity and smoothness of the velocity profile are sufficient.
Anisotropic interactions of multiple crowds lead to a non-zero-sum-game formulation, and we have developed additional criteria to guarantee the uniqueness of a Nash equilibrium for such games.
Up until now these issues were largely ignored in  the analysis of some of the more popular models \cite{Jiang_2009, Cristiani_2015},
and we showed that they can in principle lead to ill-posedness for a range of physically relevant parameter values.
However, we emphasize that our primary goal is not to criticize the prior literature, but to advocate the use of new analytic criteria -- either as a pre-processing step or as a run-time ``sanity-check'' -- in future implementations of anisotropic models.

We have discussed the numerical implementation of such models and included a number of two-crowd simulation results.  We have further proposed an approach for decreasing the computational cost of models with non-local interactions.  The intra-crowd anisotropies may arise due to the lack of radial symmetry in sensing regions \cite{Cristiani_2015},
making the computation of velocity profiles quite costly in itself. We showed how to do this efficiently wherever the crowd density is approximately linear.

We hope that our work will shift the focus of model comparison away from purely phenomenological criteria.  E.g., numerical simulations of multiple crowds can produce ``lane formation'' even if the underlying model is inconsistent or the anisotropic interactions are resolved incorrectly (using gradient descent in $u$ instead of the optimal direction field).
However, anisotropic interactions might influence the lane-orientation and are clearly important factors for the behavior of pedestrians in overlapping/interacting ``fringes'' of such lanes.

In the future, we would like to investigate several extensions both on the modeling/theoretic side and in numerical techniques.
First, the uniqueness of Nash equilibria for general (non-zero-sum) differential games is a challenging question when the players' dynamics are not separable \cite{Bressan2006}, and we hope that our criteria can be suitably generalized.  Second, a semi-Lagrangian discretization of the corresponding system of Hamilton-Jacobi-Isaacs PDEs presents another (discretized) game, and it would be useful to show that similar criteria guarantee the uniqueness of the Nash equilibrium in it. (Otherwise, the results of numerical simulations are rather hard to interpret \cite{Cacace_NZS_games}).  Finally, we would like to explore the related issues for Mean Field Games (MFG).
Up until now, the MFG-type multi-crowd models have been isotropic (and assumed the presence of stochastic perturbations in pedestrian dynamics, resulting in ``viscous'' second-order terms in PDEs) \cite{lachapelle2011mfg}.  While it is easy to write down the corresponding anisotropic/inviscid equations, we expect them to be 
similarly
affected by questions of velocity profile convexity and uniqueness of Nash equilibria.  Moreover, the standard MFG framework assumes that each agent's cost/dynamics are only affected by her own state and choice of actions + the aggregate state of all other agents \cite{lasry2006jeux, lasry2006jeux2, HuangMalhameCaines, gueant2011mfg}.  In contrast, anisotropic multi-crowd models must also account for the direct dependence on other agents' current actions, and it is not obvious whether the standard MFG framework will remain suitable.

\vspace*{2mm}
\noindent
{\bf Acknowledgements: } AV is grateful to Leighton Arnold, whose undergraduate independent study project served as a starting point for AV's work on crowd dynamics modeling.  Both authors would like to thank Emiliano Cristiani, Yanqun Jiang, and Chi-Wang Shu for very helpful discussions of their models.  The authors are also grateful to the anonymous reviewers and the editor for suggestions on improving this paper.

\vspace*{2mm}

\appendix
\numberwithin{equation}{section}
\section{Velocity Profile for Linear Density Intra-Crowd Anisotropy}
In this section, we derive a formula for the velocity profile of the model used in \cite{Cristiani_2015}, under the assumption of a linear density distribution within the sensing region.
We recall that for this model of intra-crowd anisotropy, the interactional velocity is given by:
$$
\bw_i(\x,\bu) = \int_{S(\x,\bu) \cap \Omega} \mathcal{F}(\y-\x)\rho(t,\y)d\y
$$
The original formulation, which can be found in Equation \eqref{eq:IntraCrowdSpeed} in Section \ref{ss:SpeedProfiles_Intra}, involved a cutoff value $C$ in the definition of $\mathcal{F}(\y)$. 
To simplify the analysis we shall take $C \to \infty$ and use $\mathcal{F}(\y) = \frac{-F\y}{|\y|^2}$.
Suppose that $\bu = (\cos(\theta),\sin(\theta))$, $S(\x,\bu)$ is a sector of radius $R$ and angle $\alpha$ centered in the direction $\bu$. 
Furthermore, suppose $\x$ is far enough away from $\partial\Omega$ so that the sensing region $S(\x,\bu) \subset \Omega$ and $\rho(\x)>0$ on $S(\x,\bu)$ for all choices of $\bu.$
Then fixing the time $t>0$, we may assume WLOG that our density $\rho$ is of the form $\rho(x,y) = \rho_0 + \rho_x x$, with $\rho_x >0$. 
Our formula for $\bw_i$ then becomes:

$$ \bw_i = \int_{S(\theta)} \frac{-F\y}{|\y|^2}(\rho_0+\rho_x x) d\y $$ 

which we can write in polar coordinates as:

$$ \bw_i = \int_{\theta-\alpha/2}^{\theta+\alpha/2} \int_{0}^R \frac{-F}{r^2} \left(\begin{array}{c} r\cos(\gamma) \\ r\sin(\gamma) \end{array}\right) (\rho_0 + \rho_x r\cos(\gamma))r dr d\gamma $$

After computing this integral, we get that the overall velocity $\bv = \bw_i + \bu$ is given as:

\begin{equation}\label{eq:LinearSpeed}
\bv(\theta) = \left(\begin{array}{c} \frac{-F\rho_xR^2\alpha}{4} \\ 0 \end{array}\right)
+ (1-2F\rho_0R\sin(\alpha/2)) \left(\begin{array}{c} \cos(\theta) \\ \sin(\theta) \end{array}\right)
- \frac{F\rho_xR^2\sin(\alpha)}{4} \left(\begin{array}{c} \cos(2\theta) \\ \sin(2\theta) \end{array}\right)
\end{equation}
As long as the pedestrian density is linear within the sensing region, this gives an explicit form for $\bv(\theta)$.

\label{ss:Appendix1}

\section{Convexity for Linear Density Intra-Crowd Anisotropy}
For convenience, we shall now use $C_1 = (1-2F\rho_0R\sin(\alpha/2))$, $C_2 = F\rho_xR^2\sin(\alpha)$, and $C_3 = F\rho_xR^2\alpha$. With this notation, we can write our explicit form for $\bv$ from Section \ref{ss:Appendix1} as:
\begin{equation}\label{eq:Cspeed}
\bv(\theta) = \left(\begin{array}{c} v_1(\theta) \\ v_2(\theta) \end{array}\right) 
= \left(\begin{array}{c} -\frac{C_3}{4} \\ 0 \end{array}\right)
+ C_1 \left(\begin{array}{c} \cos(\theta) \\ \sin(\theta) \end{array}\right)
- \frac{C_2}{4} \left(\begin{array}{c} \cos(2\theta) \\ \sin(2\theta) \end{array}\right).
\end{equation}

We first note that as long as $0<\alpha<\pi$, i.e. a pedestrian is only slowed down by other pedestrians in front of them, then we have $C_3>C_2>0.$

We will need to put two different conditions on our parameters so that our solutions can be logically consistent.
First, the component of the overall velocity $\bv(\theta)$ in the direction of the behavioral velocity $\bu=(\cos(\theta),\sin(\theta))$ is positive for all values of $\theta$.
We calculate this to happen when:
$$ C_1 - \frac{C_2+C_3}{4}\cos(\theta) > 0, \quad \text{ for all } \theta \in [0,2\pi].$$
Since $C_2$ and $C_3$ are always positive, this is equivalent to the condition:
\begin{equation}
C_1-\frac{C_2 + C_3}{4} > 0. 
\label{eq:origin_cond}
\end{equation}

We also need a sufficient condition for the strict convexity of the velocity profile.
As with the inter-crowd case, our velocity profile is given as a parametrized curve, so to check that it is strictly convex, we must make sure the sign of the cross product $\textbf{T} \times \textbf{N}$ does not change. 
This condition becomes:
\begin{equation}
v_1'v_2'' - v_2'v_1'' > 0, \quad \text{ for all } \theta \in [0,2\pi].
\label{eq:convex_cond}
\end{equation}
We note that Equation \eqref{eq:intercrowd_convexity} is a special case of Equation \eqref{eq:convex_cond} when $v_1(\theta)=f(\theta)\cos(\theta)$ and $v_2(\theta) = f(\theta)\sin(\theta)$.
Plugging  \ref{eq:Cspeed} into  \ref{eq:convex_cond} reduces the latter to 
\begin{equation}
C_1^2 + \frac{C_2^2}{2} + \frac{3}{2}C_1C_2\cos(\theta) > 0, \quad \text{ for all } \theta \in [0,2\pi].
\label{eq:convex_cond2}
\end{equation}

We recall that $C_2, C_3 > 0$. 
Then as long as Equation \ref{eq:origin_cond} is satisfied, we have $C_1 > 0$ as well. 
Looking at Equation \ref{eq:convex_cond2}, we see that the left hand side is minimized when $\theta = \pi$, so the sufficient condition for strict convexity becomes: 
\begin{equation}
C_1^2 + \frac{C_2^2}{2} - \frac{3}{2}C_1C_2 = \frac{(C_1 - C_2)(2C_1 - C_2)}{2} > 0.
\label{eq:convex_cond3}
\end{equation}

The case we want to avoid is when we have a non-strictly-convex velocity profile that contains the origin, i.e. when Equation \eqref{eq:origin_cond} is satisfied and Equation \eqref{eq:convex_cond3} is not.
This can only happen when:
$$ C_2 >  C_1 > \frac{C_2 + C_3}{4}. $$
Writing this in terms of the original parameters:
$$ F\rho_xR^2\sin(\alpha) > 1-2F\rho_0R\sin(\alpha/2) > F\rho_xR^2\frac{(\alpha + \sin(\alpha))}{4}.$$
These inequalities can only be satisfied if:
$ \sin(\alpha) > \frac{\alpha + \sin(\alpha)}{4} $, i.e. 
$$\alpha < 3\sin(\alpha),$$
which is true for $0 < \alpha < \approx 130.57^\circ$. 
This means that if $\alpha > \approx 130.57^\circ$ (as in \cite{Cristiani_2015}) and the velocity profile contains the origin, then the velocity profile will always be strictly convex.

\label{ss:Appendix2}

\bibliographystyle{amsplain}
\bibliography{crowds}

\end{document}